\documentclass[twoside,11pt,leqno]{article}
\usepackage{amssymb}
\usepackage{amsthm}

\usepackage{mathrsfs}
\usepackage{graphicx}
\newcommand{\zf}{\mathnormal{\mathsf{ZF}}}

\renewcommand{\max}{\mathop{\mathrm{m\acute{a}x}}}
\renewcommand{\min}{\mathop{\mathrm{m\mbox{\'{\i}}n}}}

\newtheorem{thm}{Teorema}[section]

\newtheorem{cor}[thm]{Corolario}
\newtheorem{prop}[thm]{Proposici\'on}

\theoremstyle{definition}

\newtheorem{ejem}[thm]{Ejemplo}
\newtheorem{defi}[thm]{Definici\'on}

\newcommand{\rnf}{\renewcommand{\thefootnote}{\arabic{footnote}}}
\newcommand{\thankyou}[2]{\stepcounter{footnote}\footnotetext[#1]{#2}}
\newcommand{\dem}{\noindent  {\it Demostraci\'on. }}

\title{\vspace{-3cm}
\ \hspace{-2.70in}
\\
\vspace{3cm} 
Comportamientos extra\~nos del infinito: \\
Gr\'aficas Infinitas
}

\author{\rnf  David J. Fern\'andez-Bret\'on
              \footnotemark[1]
\and  \rnf    Jes\'us A. Flores Hinostrosa
              \footnotemark[2]
\and  \rnf    V. Adrián Meza-Campa
              \footnotemark[3]
\and  \rnf    L. Gerardo N\'u\~nez Olmedo
              \footnotemark[4]
}

\date{}

\begin{document}
\maketitle


\thankyou{1}{Trabajo realizado con apoyo parcial del proyecto SIP-20221862 del IPN.}

\thankyou{2}{Trabajo realizado bajo una beca BEIFI como parte del proyecto SIP-20221862 del IPN.}

\thankyou{3}{Trabajo realizado como becario del programa Delf\'{\i}n de Investigaci\'on de Verano de 2022.}


\begin{abstract} \noindent
La combinatoria infinita (tem\'atica que, a ra\'{\i}z del trabajo de Cantor, actualmente es posible estudiar de manera completamente formal) nos presenta un interesante contraste de semejanzas y di\-fe\-ren\-cias con su an\'alogo finito. El prop\'osito de este art\'{\i}culo es presentar algunos ejemplos concretos tanto de semejanzas, como de diferencias radicales, para proporcionar cierta intuici\'on acerca del comportamiento del infinito en el \'ambito combinatorio. Nuestros ejemplos son tomados de la rama de las matem\'aticas conocida como Teor\'{\i}a de Gr\'aficas.
\end{abstract}

\noindent 


\section{Introducci\'on}
La {\em combinatoria infinita} consiste en el estudio formal de diversas estructuras infinitas, de manera an\'aloga a como se hace com\'unmente, en la combinatoria usual, con estructuras finitas. Por ejemplo, uno de los resultados b\'asicos m\'as cl\'asicos en combinatoria finita ser\'{\i}a el estudio de los coeficientes binomiales, definidos como
\begin{eqnarray*}
{n \choose k} & = & \mathrm{n\acute{u}mero\ de\ subconjuntos\ de\ }k\mathrm{\ elementos} \\
 & & \mathrm{de\ un\ conjunto\ con\ }n\mathrm{\ elementos},
\end{eqnarray*}
el cual, como todos aprendemos en nuestros cursos m\'as b\'asicos, resulta ser igual a $\displaystyle{\frac{n!}{k!(n-k)!}}$. Contar, por otro lado, la cantidad total de subconjuntos de un conjunto con $n$ elementos, implica realizar la cuenta anterior dejando que $k$ var\'{\i}e desde $0$ (el subconjunto m\'as peque\~no posible es el vac\'{\i}o, con $0$ elementos) hasta $n$ (el subconjunto m\'as grande posible ser\'{\i}a el total); obtenemos las relaciones
$$
\wp(n)=\sum_{k=0}^n {n \choose k}=2^n,
$$
en donde utilizamos el s\'{\i}mbolo $\wp(n)$ justamente para denotar el n\'umero total de subconjuntos de un conjunto con $n$ elementos. Gracias al trabajo de Cantor, hoy en d\'{\i}a podemos trabajar {\em en el para\'{\i}so}\footnote{Aqu\'{\i}, desde luego, estamos haciendo referencia a la famosa frase emitida por D. Hilbert: {\it Aus dem Paradies, das Cantor uns geschaffen, soll uns niemand vertreiben k\"onnen}, que podr\'{\i}a traducirse como {\it del para\'{\i}so que Cantor construy\'o para nosotros, nadie deber\'a expulsarnos}.} y ponderar preguntas an\'alogas respecto de las cantidades infinitas, con la posibilidad de explorar las respuestas a dichas preguntas de manera totalmente formal. Es interesante observar que esta \'area de estudio presenta un contraste entre similitudes y diferencias radicales con la combinatoria finita (tanto en el sentido de que ciertos enunciados se vuelven totalmente falsos al pasar de lo finito a lo infinito, como en el sentido de que algunos otros enunciados permanecen verdaderos, aunque la demostraci\'on co\-rres\-pon\-dien\-te resulta tener muy poco en com\'un con el caso finito). Por ejemplo, si $\kappa$ es un cardinal infinito, y $\lambda\leq\kappa$ es cualquier otro cardinal (finito o infinito), entonces se suele denotar por $[\kappa]^\lambda$ a la cardinalidad del conjunto cuyos elementos son los subconjuntos de cardinalidad $\lambda$ de un conjunto dado de tama\~no $\kappa$ (en otras palabras, este objeto es el an\'alogo infinito del coeficiente binomial). En ciertos casos, resulta bastante complicado determinar con exactitud el valor de $[\kappa]^\lambda$: si bien es cierto que $[\kappa]^0=1$, y $[\kappa]^\lambda=\kappa$ siempre que $\lambda\neq 0$ sea finito, el valor de $[\kappa]^\lambda$ para $\lambda<\kappa$ infinito es, en cierto sentido, un completo interrogante, pues el valor exacto var\'{\i}a dependiendo del tipo de cardinal que sea $\kappa$, y en muchos casos este valor es incluso --demostrablemente-- imposible de determinar utilizando \'unicamente los axiomas usuales de la teor\'{\i}a de conjuntos (para valores de $\lambda>\kappa$, sin embargo, es evidente que $[\kappa]^\lambda=0$, y por otra parte tambi\'en es cierto que $[\kappa]^\kappa=2^\kappa$). Desde luego que esta situaci\'on contrasta significativamente con lo que sucede en el caso finito. Sin embargo, sigue siendo cierto que $\displaystyle{\sum_{\lambda\leq\kappa}[\kappa]^\lambda}=2^\kappa$, aunque en este caso la demostraci\'on ya presenta discrepancias importantes con la correspondiente demostraci\'on del caso finito.

El prop\'osito de este art\'{\i}culo divulgativo es presentar al lector con varios ejemplos, est\'eticamente muy satisfactorios, de estos contrastes entre la combinatoria finita y la infinita. Como terreno de juego para dicha actividad recreativa, hemos elegido el f\'ertil campo de la Teor\'{\i}a de Gr\'aficas, que constituye una de las \'areas fundamentales de la combinatoria en el caso finito --y que tambi\'en en su momento inspir\'o varias investigaciones fundamentales en el caso infinito--. Para esto, ser\'{\i}a de gran ayuda que el lector se encuentre ya ligeramente familiarizado con las vicisitudes de la aritm\'etica de cardinales infinitos, aunque en la secci\'on 2 incluimos un breve sumario de las definiciones y resultados relevantes (dichos resultados s\'olo se mencionan, y no se demuestran, pero el lector con madurez matem\'atica deber\'{\i}a poder seguir el paso del resto del art\'{\i}culo si simplemente toma el contenido de esta secci\'on como art\'{\i}culo de fe). En la secci\'on 3 incluimos algunos resultados ``b\'asicos'', en el sentido de que pueden enunciarse con s\'olo las definiciones b\'asicas de teor\'{\i}a de gr\'aficas, sin involucrar a\'un demasiada teor\'{\i}a. La secci\'on 4 versa sobre los \'arboles y bosques (definidos exactamente igual a como se hace en el caso finito). La secci\'on 5 introduce las coloraciones por v\'ertices y el n\'umero crom\'atico. En la secci\'on 6 exploramos cuestiones que tienen que ver con emparejamientos de v\'ertices, as\'{\i} como la relaci\'on de esto con algunos de los temas tratados anteriormente. Esta \'ultima secci\'on ilustra con abundancia un fen\'omeno bastante interesante que sucede con frecuencia al movernos del terreno finito al infinito: numerosos teoremas ``siguen siendo ciertos'' en el contexto infinito, siempre y cuando se encuentre la manera correcta de enunciarlos en este contexto m\'as general. Esta secci\'on es tambi\'en aquella en la que m\'as se han realizado investigaciones de alto grado de dificultad, y se encuentran varios teoremas cuyas demostraciones constituyen la mayor parte de art\'{\i}culos completos de investigaci\'on; muchos de estos resultados \'unicamente se mencionan y se ponen en contexto, sin incluir sus demostraciones. Finalmente, se incluye al final del art\'{\i}culo un ap\'endice que contiene algunos contraejemplos en el contexto donde no se asume el llamado {\em Axioma de Elecci\'on}, mismo que juega un papel importante en varias de las demostraciones contenidas en las secciones previas; de esta forma, se explica con cierto detalle la manera como este axioma es esencial a las demostraciones mencionadas.

Intentamos presentar las demostraciones pertinentes con un nivel de detalle que permita a las ideas principales quedar expuestas, sin permitir que los excesivos detalles mec\'anicos afecten ni la legibilidad, ni la visibilidad del flujo principal de ideas. En particular, los lectores menos experimentados que quieran realmente comprender a profundidad lo aqu\'{\i} expuesto, deber\'an tratar cada afirmaci\'on no demostrada como si fuera un ejercicio. La frase ``es f\'acil ver'', as\'{\i} como cualquier otra frase a efecto similar, funcionan como se\~nalamientos de la pre\-sen\-cia de alg\'un ejercicio de este tipo (ejercicios que, si bien no siempre son precisamente f\'aciles, s\'{\i} son, en la mayor\'{\i}a de los casos, rutinarios, en el sentido de que constituyen alguna verificaci\'on mec\'anica de la veracidad de alguna proposici\'on sencilla, sin necesidad de involucrar ideas novedosas de ning\'un tipo).

\section{Aritm\'etica cardinal infinita}\label{sect:aritmetica}

La definici\'on clave, descubierta por Cantor (justo esa idea genial gracias a la cual hoy en d\'{\i}a podemos estudiar formalmente el infinito), es la siguiente: es posible comparar las cantidades de elementos de dos conjuntos sin necesidad de contar a los elementos de cada uno de estos conjuntos por separado; antes bien, basta emparejar a los elementos de un conjunto con los del otro, sin permitir que nadie se quede sin emparejar y que nadie tenga m\'as de una pareja.

\begin{defi}\label{equipotencia}
Dados dos conjuntos $X,Y$, decimos que
\begin{enumerate}
\item $X$ es {\em equipotente a} $Y$, simbolizado $X\approx Y$ o $|X|=|Y|$, si existe una funci\'on biyectiva entre $X$ y $Y$.
\item $X$ {\em tiene a lo m\'as tantos elementos como $Y$}, simbolizado $X\preceq Y$ o $|X|\leq |Y|$, si existe una funci\'on inyectiva $f:X\longrightarrow Y$.
\end{enumerate}
\end{defi}

Deber\'{\i}a ser bastante claro para el lector que la relaci\'on de equipotencia $\approx$ se comporta como una relaci\'on de equivalencia (es reflexiva, sim\'etrica y transitiva); el \'unico detalle es que no est\'a definida sobre un conjunto, sino sobre la clase (propia) de todos los conjuntos, pero en este art\'{\i}culo nos abstendremos de discutir esas sutilezas fundacionales de la teor\'{\i}a axiom\'atica de conjuntos. Es por ello que nos abstenemos de asignarle una denotaci\'on concreta al s\'{\i}mbolo $|X|$ (que, intuitivamente, deber\'{\i}a ser ``la cardinalidad de $X$'', la clase de equivalencia de $X$ m\'odulo la relaci\'on de equipotencia), y jam\'as lo utilizaremos por s\'{\i} solo; \'unicamente utilizaremos los s\'{\i}mbolos $|X|=|Y|$ y/o $|X|\leq|Y|$ como abreviaciones de los enunciados correspondientes, tal y como vienen estipulados en la Definici\'on~\ref{equipotencia}.

Por otra parte, la relaci\'on $\preceq$ es antisim\'etrica y transitiva; es decir, se comporta como una relaci\'on de preorden --con lo cual, el comportamiento es como el de una relaci\'on de orden parcial en las clases de equivalencia m\'odulo la relaci\'on que hace a $X$ y $Y$ equivalentes si y s\'olo si $X\preceq Y$ y $Y\preceq X$--. Para comprender esta relaci\'on de equivalencia, notemos que, dados dos conjuntos $X,Y$, se tiene que $X\preceq Y$ y $Y\preceq X$ si y s\'olo si $X\approx Y$; una de estas implicaciones es obvia, y la otra dista mucho de ser trivial, al grado de que es un teorema que tiene nombre: el teorema de Cantor--Bernstein.

La aritm\'etica de n\'umeros cardinales infinitos se define por analog\'{\i}a con el caso de los cardinales finitos. Dados dos cardinales $\kappa$ y $\lambda$, definimos a $\kappa+\lambda$ como la cardinalidad de la uni\'on disjunta entre un conjunto de cardinalidad $\kappa$, y otro de cardinalidad $\lambda$ (como ya dijimos arriba, al escribir n\'umeros cardinales como si fueran objetos, uno debe en realidad pensar en estos enunciados como abreviaturas: as\'{\i}, nuestra ``definici\'on'' de suma en realidad debe concebirse como la estipulaci\'on de que, siempre que escribamos $|Z|=|X|+|Y|$, lo que en realidad estamos diciendo es que $Z$ es equipotente a la uni\'on disjunta de $X$ y $Y$). Similarmente, definimos $\kappa\lambda$ como la cardinalidad del producto cruz $X\times Y$, con $|X|=\kappa$ y $|Y|=\lambda$; y finalmente definimos $\kappa^\lambda$ como la cardinalidad del conjunto $X^Y=\{f\big|f:Y\longrightarrow X\}$, en donde $|X|=\kappa$ y $|Y|=\lambda$. El lector podr\'a sin duda verificar, sin demasiada dificultad, que las tres definiciones coin\-ci\-den con las definiciones usuales de la aritm\'etica elemental en el caso de los cardinales finitos. El teorema fundamental que caracteriza a la suma y el producto de cardinales, es el siguiente.

\begin{thm}
Sean $\kappa,\lambda$ dos n\'umeros cardinales, en donde al menos uno de los dos es infinito. Entonces, se cumple lo siguiente:
\begin{enumerate}
\item $\kappa+\lambda=\max\{\kappa,\lambda\}$,
\item si $\kappa,\lambda$ son ambos no cero, entonces $\kappa\lambda=\max\{\kappa,\lambda\}$.
\end{enumerate}
\end{thm}

Desde luego, si tanto $\kappa$ como $\lambda$ son finitos en el teorema anterior, entonces el resultado tanto de la suma como del producto es, como hemos remarcado anteriormente, ya conocido: coincide con la suma y el producto de n\'umeros naturales. Similarmente, si alguno de los dos n\'umeros, $\kappa$ o $\lambda$, es igual a cero, entonces $\kappa\lambda=0$ (independientemente de que el otro n\'umero cardinal sea finito o infinito).

Tambi\'en es posible realizar tanto sumatorias como productos de una familia de cardinales. Esto es, si $\{X_\alpha\big|\alpha\in\Lambda\}$ es una familia indexada de conjuntos, con conjunto de \'{\i}ndices $\Lambda$, entonces se define $\sum_{\alpha\in\Lambda}|X_\alpha|$ como la cardinalidad de la uni\'on disjunta de los conjuntos $X_\alpha$; formalmente,
$$
\sum_{\alpha\in\Lambda}|X_\alpha|=\bigg|\bigcup_{\alpha\in\Lambda}(X_\alpha\times\{\alpha\})\bigg|.
$$
Similarmente, $\prod_{\alpha\in\Lambda}|X_\alpha|=\big|\prod_{\alpha\in\Lambda}X_\alpha\big|$. Es posible demostrar que, siempre que los $|X_\alpha|$ sean todos ellos distintos de cero, se tiene que
$$
\sum_{\alpha\in\Lambda}|X_\alpha|=\sup(\{|X_\alpha|\big|\alpha\in\Lambda\}\cup\{|\Lambda|\});
$$
por otra parte, el llamado {\em lema de K\"onig} nos garantiza que, si $\{X_\alpha\big|\alpha\in\Lambda\}$ y $\{Y_\alpha\big|\alpha\in\Lambda\}$ son sendas familias indexadas de conjuntos, tales que $|X_\alpha|<|Y_\alpha|$ para todo $\alpha\in\Lambda$, entonces
$$
\sum_{\alpha\in\Lambda}|X_\alpha|<\prod_{\alpha\in\Lambda}|Y_\alpha|.
$$

Se utilizar\'an tambi\'en, en mucha menor medida, los n\'umeros ordinales. Estos fueron originalmente descubiertos por Cantor, aunque la formalizaci\'on m\'as utilizada actualmente es la de von Neumann, en donde cada n\'umero ordinal se identifica con el conjunto de aquellos n\'umeros ordinales menores a \'el (por ejemplo, $0=\varnothing$, $1=\{0\}$ y, en general, $n=\{0,\ldots,n-1\}$; la formalizaci\'on de von Neumann permite continuar este conteo hasta el \'ambito transfinito). De esta forma, los e\-le\-men\-tos de un n\'umero ordinal se encuentran bien ordenados mediante la relaci\'on de orden parcial $\subseteq$, en donde la relaci\'on de pertenencia $\in$ cons\-ti\-tu\-ye el orden estricto correspondiente. Cada conjunto bien ordenado es isomorfo a un \'unico n\'umero ordinal; de esta forma, los n\'umeros ordinales nos proporcionan una colecci\'on que contiene a un ``representante can\'onico'' de cada ``clase de equivalencia'' de conjuntos bien ordenados m\'odulo la relaci\'on de isomorfismo. Por otra parte, la clase misma de los n\'umeros ordinales est\'a tambi\'en bien ordenada, en el sentido de que cada conjunto no vac\'{\i}o de ordinales tiene un elemento $\subseteq$-m\'{\i}nimo; esto permite realizar demostraciones por inducci\'on transfinita y definir objetos por recursi\'on transfinita, es decir, si realizamos alg\'un procedimiento en un n\'umero ordinal $\alpha$, utilizando como suposici\'on que el procedimiento ya ha sido realizado con todos los ordinales $\beta<\alpha$, entonces este pro\-ce\-di\-mien\-to (ya sea una definici\'on, o una demostraci\'on) est\'a bien definido. Adem\'as, si suponemos el axioma de elecci\'on, entonces todo conjunto puede ser bien ordenado; luego, todo conjunto es equipotente a alg\'un n\'umero ordinal, y por lo tanto existe un m\'{\i}nimo n\'umero ordinal equipotente a nuestro conjunto. Como consecuencia de esto, siempre que exista alg\'un conjunto con cierta propiedad, es v\'alido tomar un tal conjunto de cardinalidad m\'{\i}nima posible, lo cual se utilizar\'a con bastante frecuencia en diversas definiciones a lo largo del presente art\'{\i}culo.

\section{Resultados b\'asicos}

Sin suponer teor\'{\i}a previa alguna, comenzamos desde la definici\'on misma de gr\'afica, que proporcionamos a continuaci\'on. El lector con algo de experiencia notar\'a que, a lo largo de este art\'{\i}culo, nuestras gr\'aficas son todas ellas simples y no-dirigidas.

\begin{defi}
Una {\em gr\'afica} es una pareja ordenada $G=(V,E)$, en donde $V$ es cualquier conjunto no vac\'io, y\footnote{En teor\'{\i}a de conjuntos, es est\'andar utilizar la notaci\'on $[X]^\kappa$ para referirse al conjunto de subconjuntos de $X$ que tienen cardinalidad $\kappa$. As\'{\i}, la presente definici\'on estipula que cada elemento de $E$ es un conjunto que contiene a exactamente dos elementos de $V$. De esta forma, identificamos a una arista --que, en el presente art\'{\i}culo, nunca es un bucle-- con el conjunto de sus dos extremos.} $E\subseteq[V]^2$. A los elementos de $V$ se les llama {\em v\'ertices}, y a los de $E$ se les llama {\em aristas}.
\end{defi}

Notoriamente, nuestra definici\'on de gr\'afica (simple y no-dirigida) coincide con la que usualmente se utiliza en el caso finito, con la \'unica excepci\'on de que nos abstenemos de insistir en este \'ultimo punto --es decir, simplemente no requerimos que el conjunto de v\'ertices debe de ser finito--.

Un concepto est\'andar en teor\'{\i}a de gr\'aficas es el del {\em grado} de un v\'ertice. En el caso que nos ocupa (donde no hay bucles ni aristas m\'ultiples), es f\'acil definir el grado de un v\'ertice como la cantidad, ya sea de v\'ertices adyacentes al dado, o bien de aristas incidentes al v\'ertice en cuesti\'on. En s\'{\i}mbolos, si $G=(V,E)$ es una gr\'afica, y $v\in V$, entonces definimos
$$
d_G(v)=|\{u\in V\big|\{u,v\}\in E\}|=|\{e\in E\big|v\in e\}|.
$$

Uno de los resultados m\'as elementales en teor\'{\i}a de gr\'aficas finitas es el que afirma que, al sumar los grados de todos los v\'ertices de una gr\'afica, uno obtiene el doble del n\'umero de aristas que contiene la gr\'afica. En el caso finito, no es dif\'{\i}cil formalizar el argumento usual, el cual se basa en la observaci\'on de que, al sumar los grados de cada v\'ertice, se est\'a al mismo tiempo realizando un conteo de aristas, con la salvedad de que cada arista se cuenta dos veces (una por cada uno de sus extremos). En el caso infinito, sin embargo, es necesario prestar un poco m\'as de atenci\'on al detalle, como a continuaci\'on veremos.

\begin{thm}\label{teor:2e}
Si $G=(V,E)$ es cualquier gr\'afica, entonces
$$
\sum_{v \in V} d_G(v) = 2|E|.
$$
\end{thm}

\dem
Dadas las definiciones proporcionadas en la secci\'on~\ref{sect:aritmetica}, lo que en realidad necesitamos demostrar es que, si $X_v$ es un conjunto de cardinalidad $d_G(v)$ para cada $v\in V$, de modo que los $X_v$ sean disjuntos a pares, entonces $\bigcup_{v\in V}X_v$ es equipotente a $E\times\{0,1\}$ (de acuerdo con la construcci\'on de von Neumann, usualmente se toma a $\{0,1\}$ como el conjunto paradigm\'atico de cardinalidad 2). Lo primero que se viene a la mente al momento de elegir a $X_v$ es dejar que sea o bien el conjunto de vecinos de $v$, o bien el conjunto de aristas que inciden en $v$; sin embargo, con cualquiera de estas definiciones nos enfrentar\'{\i}amos al problema de que los $X_v$ no ser\'{\i}an disjuntos por pares. De esta forma, unos momentos de reflexi\'on nos llevan a proponer los conjuntos
$$
X_v=\{(v,u)\big|\{v,u\}\in E\}.
$$
La clave es que $(v,u)$ es un par ordenado, y no simplemente un conjunto con dos elementos (de modo que $(v,u)$ y $(u,v)$ son dos pares ordenados distintos, ambos asociados a la arista $\{u,v\}$). Esto nos da la tranquilidad de saber que $X_v\cap X_u=\varnothing$ siempre que $v\neq u$.

Utilizando el axioma de elecci\'on, tomamos una funci\'on de elecci\'on $f:E\longrightarrow V$; esto es, una funci\'on tal que $f(e)\in e$ para todo $e\in E$. En otras palabras, $f$ es un objeto que representa el haber elegido uno de los dos extremos de cada una de las aristas de la gr\'afica $G$. As\'{\i}, definimos la funci\'on $\varphi:\bigcup_{v\in V}X_v\longrightarrow E\times\{0,1\}$ mediante la f\'ormula
$$
\varphi(u,v)=\ \left\{\begin{array}{l}
(\{u,v\},0),\ \ \ \ \ \mathrm{\ si\ }f(\{u,v\})=u; \\
(\{u,v\},1),\ \ \ \ \ \mathrm{\ si\ }f(\{u,v\})=v.
\end{array}\right.
$$
No es precisamente trivial, pero s\'{\i} rutinario, demostrar que la funci\'on $\varphi$ definida arriba es una biyecci\'on.

\hfill$\Box$

En la demostraci\'on anterior, la principal dificultad t\'ecnica es que, para cada arista $\{u,v\}$, tenemos dos copias de dicha arista (co\-rres\-pon\-dien\-tes a los pares ordenados $(u,v)$ y $(v,u)$) y debemos decidir, al definir nuestra biyecci\'on $\varphi$, cu\'al de estas dos copias se manda a la pareja cuya segunda entrada es $0$, y cu\'al a la pareja con segunda entrada $1$. De alguna manera, estamos eligiendo un elemento dentro del conjunto $\{(u,v),(v,u)\}$ (estamos eligiendo, por ejemplo, cu\'al de los dos debe ser enviado a la pareja con segunda entrada igual a $1$). De modo que el axioma de elecci\'on juega un papel crucial en la demostraci\'on previa; el ap\'endice al final de este art\'{\i}culo muestra que este axioma es, de hecho, indispensable, pues mostramos un ejemplo concreto (Ejemplo~\ref{ex:secondfraenkel}), en la teor\'{\i}a axiom\'atica de Zermelo--Fraenkel sin el axioma de elecci\'on, en el cual el Teorema~\ref{teor:2e} es falso.

\begin{figure}\label{fig:camino-infinito}
\begin{center}
\includegraphics[scale=.6]{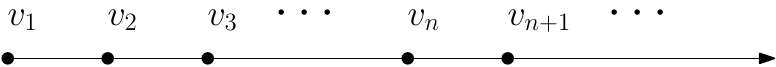}
\end{center}
\caption{El camino infinito en una direcci\'on.}
\label{fig:cam-infinito}
\end{figure}

A continuaci\'on analizamos dos resultados que muestran un contraste interesante entre la combinatoria finita y la infinita. En primer lugar, tras haber demostrado el Teorema~\ref{teor:2e}, en el caso de las gr\'aficas finitas es inmediato verificar, como corolario, que la cantidad de v\'ertices de grado impar debe de ser un n\'umero par (ya que la suma de todos estos grados, en total, debe de dar $2|E(G)|$, que es un n\'umero par, y la suma de una cantidad impar de n\'umeros impares siempre resulta en un n\'umero impar). En el caso infinito, esto no necesariamente se cumple, como lo ilustra el ejemplo de un camino infinito en una direcci\'on. En este caso, nuestra gr\'afica $G=(V,E)$ cuenta con una cantidad numerable de v\'ertices, digamos que $V=\{v_1,v_2,\dots,v_n,\ldots\}$, y hacemos $E= \{\{v_n,v_{n+1}\}\big|n\in\mathbb N\}$. As\'{\i} obtenemos una gr\'afica (re\-pre\-sen\-ta\-da gr\'aficamente en la Figura~\ref{fig:cam-infinito}) con exactamente un v\'ertice de grado 1, y los dem\'as de grado 2.

En segundo lugar, tenemos otro resultado sencillo (que usualmente suele ser un ejercicio en los textos introductorios de teor\'{\i}a de gr\'aficas): para cualquier gráfica (finita) con al menos dos vértices, siempre hay dos vértices con el mismo grado. La demostraci\'on de este resultado finito reposa en el principio finito de la pichonera, de la manera siguiente: si suponemos que nuestra gr\'afica tiene $n$ v\'ertices, entonces (dado que \'unicamente estamos considerando gr\'aficas simples) cada uno de esos v\'ertices debe tener un grado (igual al n\'umero de sus vecinos) entre $0$ y $n-1$. Si no hubiera dos v\'ertices con el mismo grado, entonces los grados de cada uno de los v\'ertices ser\'{\i}an exactamente los n\'umeros $0,1,\ldots,n-1$, sin repeticiones; esto ya es una contradicci\'on, pues no pueden existir simult\'aneamente un v\'ertice de grado $0$ (ning\'un vecino) y uno de grado $n-1$ (todos los dem\'as v\'ertices como vecinos). La demostraci\'on que acabamos de proporcionar est\'a profundamente enraizada en la finitud de la gr\'afica en cuesti\'on, por lo que no resulta sorprendente que, en el caso infinito, el an\'alogo de este resultado sea falso. Para construir un contraejemplo, podemos tomar una gr\'afica con conjunto de v\'ertices $V = \{v_n\big|n\in\mathbb N\}$, en donde el conjunto de aristas se define mediante $E = \{\{v_n,v_{n+k}\}\big|0< k \leq 2^n \}$. N\'otese que, para cada $n$, el v\'ertice $v_n$ \'unicamente puede ser adyacente a v\'ertices de la forma $v_j$ con $1\leq j\leq n+2^n$; as\'{\i}, tenemos que $2^n\leq d(v_n)\leq 2^n+n<2^{n+1}\leq d(v_{n+1})$. Por lo tanto, no es dif\'{\i}cil verificar que la asignaci\'on $n\longmapsto d(v_n)$ es inyectiva; en otras palabras, no hay dos v\'ertices de esta gr\'afica que tengan el mismo grado.

Recordemos que, en el caso de una gr\'afica finita $G$, se definen los n\'umeros $\delta(G)$ y $\Delta(G)$ como los grados m\'{\i}nimo y m\'aximo, respectivamente, de alg\'un v\'ertice de la gr\'afica. En el caso infinito, el conjunto de estos grados podr\'{\i}a ser no-acotado; o inclusive podr\'{\i}a darse el caso de que algunos de los grados sean cardinales infinitos. Esto nos fuerza a modificar la definici\'on del grado m\'aximo (la del grado m\'{\i}nimo, por otra parte, se puede mantener, gracias al resultado (de hecho, equivalente al axioma de elecci\'on) que afirma que la clase de todos los n\'umeros cardinales est\'a bien ordenada).

\begin{defi}
Sea $G=(V,E)$ una gr\'afica (finita o infinita). Se definen los siguientes dos par\'ametros:
\begin{enumerate}
\item $\delta(G)=\min\{d_G(v)\big|v\in V\}$,
\item $\Delta(G)=\sup\{d_G(v)\big|v\in V\}$.
\end{enumerate}
\end{defi}

A continuaci\'on enunciamos un resultado junto con uno de sus corolarios. El resultado afirma, en el caso infinito, exactamente lo mismo que en el caso finito; por su parte, el corolario presenta un fen\'omeno completamente nuevo que es exclusivo del dominio infinito.

\begin{thm}
Si $G=(V,E)$ es cualquier gr\'afica, entonces $|V|\delta (G)  \leq 2|E| \leq |V| \Delta(G)$.
\end{thm}

\dem
Dado que $\delta(G)\leq d_G(v)\leq \Delta(G)$ para todo $v\in V$, entonces tenemos que
$$
|V| \delta (G)=\sum_{v\in V}\delta(G)\leq \sum_{v \in V} d_G(v)\leq \sum_{i \in V} \Delta(G)=|V|\Delta (G),
$$
en donde todas las (des)igualdades son inmediatas en el caso finito, y en el caso infinito, la primera y \'ultima igualdad se siguen de
$$
\sum_{v\in V}\delta(G)=\sup\{|V|,\delta(G)\}=\max\{|V|,\delta(G)\}=|V|\delta(G),
$$
y
$$
\sum_{v\in V}\Delta(G)=\sup\{|V|,\Delta(G)\}=\max\{|V|,\Delta(G)\}=|V|\Delta(G),
$$
respectivamente.

\hfill$\Box$

\begin{cor}\label{inf-v=e}
Sea $G=(V,E)$ una gr\'afica infinita tal que $\delta(G)>0$. Entonces, $|V|=|E|$.
\end{cor}

\dem
Como \'unicamente consideramos gr\'aficas simples, $\Delta(G)\leq|V|$. Entonces (dado que $\delta(G)$ es positivo, y que $|V|$ es infinito) tenemos que:
\begin{eqnarray*}
|V| & = & \max\{|V|,\delta(G)\}= |V| \delta (G)  \leq 2|E| \\
 & & \leq |V| \Delta(G)=\max\{|V|,\Delta(G)\} = |V|,
\end{eqnarray*}
lo cual implica que $|V| = |E|$.

\hfill$\Box$

Desde luego que, en general, el corolario anterior dif\'{\i}cilmente se cumplir\'a para gr\'aficas finitas (excepto en casos muy especiales, como por ejemplo las gr\'aficas $C_n$). Tambi\'en es claro que la hip\'otesis $\delta(G)>0$ es esencial, pues de lo contrario, es posible a\~nadir a una gr\'afica tantos v\'ertices aislados como sea necesario para forzar a que el conjunto de v\'ertices tenga una cardinalidad tan grande como se desee, sin alterar la cardinalidad del conjunto de aristas.

\section{\'Arboles}\label{sect:arboles}



Dada una gr\'afica $G=(V,E)$ y dos v\'ertices $u,v\in V$, un {\bf $(u,v)$-camino} es una sucesi\'on finita de v\'ertices $(v_0,\ldots,v_n)$ tal que $v_0=u$, $v_n=v$ y $\{v_i,v_{i+1}\}\in E$ para cada $i\in\{0,\ldots,n-1\}$. En este caso, decimos que $n$ es la {\bf longitud} del camino. Decimos que dos v\'ertices $u,v\in V$ est\'an {\bf conectados} si y s\'olo si existe un $(u,v)$-camino en $G$. Todas estas definiciones son exactamente iguales a las utilizadas en el caso de gr\'aficas finitas; al igual que en ese caso, la relaci\'on de estar conectado constituye una relaci\'on de equivalencia entre los v\'ertices. A las subgr\'aficas de $G$ generadas por cada una de las clases de equivalencia se les llama {\bf componentes conexas}, y decimos que $G$ es {\bf conexa} si y s\'olo si tiene exactamente una componente conexa. A un $(u,u)$ camino de longitud no-cero le llamamos un {\bf ciclo} de $G$, y decimos que $G$ es {\bf ac\'{\i}clica} si no admite ciclos.

\begin{defi}
A una gr\'afica ac\'{\i}clica se le llama un {\bf bosque}. Un bosque conexo se denomina un {\bf \'arbol}.
\end{defi}

La terminolog\'{\i}a forestal, adem\'as de ser sugestiva, est\'a plenamente justificada, ya que con estas definiciones un bosque constituye la uni\'on de \'arboles (cada una de las componentes conexas). Notacionalmente, es com\'un utilizar la letra $T$ para denotar \'arboles, debido a la palabra ``tree'' que significa \'arbol en ingl\'es.

A continuaci\'on enunciamos algunos hechos b\'asicos acerca de \'arboles; todas las afirmaciones del presente p\'arrafo son resultados cl\'asicos en el caso de los \'arboles finitos, y se demuestran de manera exactamente id\'entica en el caso de los \'arboles infinitos. Dado un \'arbol $T$, cualesquiera dos v\'ertices est\'an conectados mediante un \'unico camino (la existencia del camino se debe a la conexidad, mientras que la unicidad del mismo se puede demostrar con facilidad utilizando el hecho de que $T$ no contiene ciclos). En particular, dada una gr\'afica conexa $G$, tenemos que $G$ es un \'arbol si y s\'olo si todas sus aristas son de corte (una arista $e$ en una gr\'afica $G$ se dice {\em de corte} si ambos extremos de $e$ pertenecen a componentes conexas distintas de $G-e$, la gr\'afica que resulta de eliminar al elemento $e$ del conjunto de aristas de $G$ manteniendo el mismo conjunto de v\'ertices, y es f\'acil demostrar --el mismo argumento funciona tanto en el caso finito como en el infinito-- que una arista es de corte si y s\'olo si no est\'a contenida en ning\'un ciclo de $G$). Consideraremos ahora los v\'ertices de corte: un v\'ertice en una gr\'afica $G=(V,E)$ es {\em de corte} si existe una partici\'on $E=E_1\cup E_2$ de las aristas de $G$ tal que $v$ es el \'unico v\'ertice com\'un a las subgr\'aficas generadas por $E_1$ y $E_2$. Esto es, en el caso finito, equivalente a decir que $G-v$ tiene (estrictamente) m\'as componentes conexas que $G$; no as\'{\i} en el caso infinito, como se puede ver, por ejemplo, tomando una infinidad de gr\'aficas finitas, $G_n$ para $n\in\mathbb N$, y dejando que $G$ sea la uni\'on disjunta de las $G_n$: cualquier v\'ertice de corte de alguna $G_n$ sigue siendo v\'ertice de corte de $G$, aunque la cardinalidad del conjunto de componentes conexas permanece constante (con valor $\aleph_0$). Sin embargo, si $G$ es una gr\'afica conexa, entonces un v\'ertice $v$ de $G$ es de corte si y s\'olo si $G-v$ es disconexa. En el caso particular de un \'arbol $T$, un v\'ertice $v$ de $T$ es de corte si y s\'olo si $d(v)>1$.

En contraste, un resultado cl\'asico acerca de \'arboles finitos $T=(V,E)$ es que siempre satisfacen $|E|-1=|V|$. Esto es, en cierto sentido, trivialmente cierto para el caso de \'arboles infinitos: dado un \'arbol infinito $T=(V,E)$, al no haber v\'ertices aislados (pues de lo contrario $T$ no ser\'{\i}a conexo), se satisface la igualdad $|E|=|V|=|V|+1$, por el corolario~\ref{inf-v=e}. Sin embargo, un enunciado m\'as fuerte, m\'as informativo, y que presumiblemente refleja una m\'as adecuada generalizaci\'on de la igualdad $|E|-1=|V|$, es el siguiente.

\begin{thm}\label{teor:arboleigualvmasuno}
Sea $T=(V,E)$ un \'arbol (finito o infinito), y sea $v_0\in V$ arbitrario. Entonces existe una biyecci\'on $\varphi:V\setminus\{v_0\}\longrightarrow E$ tal que, para todo $v\in V\setminus\{v_0\}$, $\varphi(v)$ es una arista incidente en $v$.
\end{thm}

\dem
Para cada $v\in V\setminus\{v_0\}$, sea $(v_0,\ldots,v_n)$, con $v=v_n$, el \'unico $(v_0,v)$-camino, y denotemos $\varphi(v)=\{v_{n-1},v_n\}$ (esto tiene sentido ya que, si $v\neq v_0$, entonces $n>0$). Es rutinario verificar que la funci\'on $\varphi:V\setminus\{v_0\}\longrightarrow E$ satisface lo afirmado en el enunciado del teorema.

\hfill$\Box$

En el caso de los \'arboles finitos, el resultado $|E|=|V|-1$ permite mostrar que todo \'arbol finito con al menos dos v\'ertices contiene al menos dos v\'ertices de grado 1 (pues de lo contrario, tendr\'{\i}amos $d_G(v)\geq 2$ para todos, excepto posiblemente uno, de los $v\in V$, y por lo tanto $2|V|-2=2|E|=\sum_{v\in V}d_G(v)\geq 2(|V|-1)+1=2|V|-1$, una contradicci\'on). Esta \'ultima afirmaci\'on resulta falsa en el caso de los \'arboles infinitos, como se puede apreciar considerando el ejemplo del camino infinito en una direcci\'on, ilustrado en la Figura~\ref{fig:cam-infinito}. Como consecuencia de ello, tambi\'en es falso, en el \'ambito infinito, el enunciado (que s\'{\i} se cumple en el caso finito) que afirma que toda gr\'afica conexa con m\'as de un v\'ertice tiene cuando menos dos v\'ertices que no son de corte (nuevamente, el camino infinito en una direcci\'on funciona como contraejemplo).

Ahora pasamos a hablar de \'arboles generadores en una gr\'afica. Dada una gr\'afica $G=(V,E)$, decimos que $T$ es un {\em \'arbol generador} si es una subgr\'afica generadora de $G$ (es decir, el conjunto de v\'ertices de $T$ es $V$ mismo, y el conjunto de aristas de $T$ es un subconjunto de $E$) que adem\'as es un \'arbol. En particular, toda gr\'afica que admite un \'arbol generador debe, al contener una subgr\'afica inducida\footnote{Recuerde que, dado un subconjunto de v\'ertices $X$ de una gr\'afica $G$, la {\bf subgr\'afica inducida por $X$} es la subgr\'afica de $G$ cuyo conjunto de v\'ertices es $X$ y cuyo conjunto de aristas consta de absolutamente todas las aristas de $G$ con ambos extremos en $X$.} conexa, ser ella misma conexa. Resulta ser que la  implicaci\'on rec\'{\i}proca tambi\'en se cumple, y toda gr\'afica conexa admite un \'arbol generador. En el caso finito, esto se puede demostrar de manera m\'as o menos constructiva --una opci\'on es partir de $G$, y recursivamente eliminar aristas, una por una, de tal forma que la gr\'afica resultante siga siendo conexa (en el momento en que deja de ser posible eliminar una arista m\'as, significa que la subgr\'afica que hemos obtenido hasta ese momento es un \'arbol); la otra opci\'on es partir de la subgr\'afica generadora de $G$ sin aristas, y progresivamente ir a\~nadiendo aristas sin a\~nadir ciclos (en el momento en que ya no es posible a\~nadir una arista m\'as, significa que la subgr\'afica obtenida hasta ahora, que es un \'arbol, ya es generadora). En el caso de gr\'aficas infinitas, sin embargo, pese a que el resultado se mantiene verdadero, la demostraci\'on es notoriamente distinta --y significativamente no-constructiva, como veremos a continuaci\'on--.

\begin{thm}\label{teor:conexaimplicaarbolgen}
Toda gr\'afica conexa admite un \'arbol generador.
\end{thm}

\dem
Dada la gr\'afica conexa $G$, consid\'erese el conjunto
$$
\mathbb P=\left\{H\big|H\mathrm{\ es\ una\ subgr\acute{a}fica\ de\ }G\mathrm{\ y\ es\ un\ \acute{a}rbol}\right\},
$$
parcialmente ordenado mediante la relaci\'on $H\leq H'$ si y s\'olo si $H'$ es una subgr\'afica de $H$. N\'otese que $\mathbb P$ es no vac\'{\i}o, ya que por lo menos contiene, para cualquier v\'ertice de $G$, a la subgr\'afica (sin aristas) inducida por dicho v\'ertice. Adem\'as, si suponemos que $\mathscr C\subseteq\mathbb P$ es una cadena, entonces la uni\'on $\tilde{H}$ de los elementos de $\mathscr C$ es tambi\'en un elemento de $\mathbb P$, como argumentaremos a continuaci\'on. Si $u,v$ son v\'ertices de $\tilde{H}$, entonces $u,v$ son v\'ertices de alg\'un $H\in\mathscr C$ (inicialmente, $u$ y $v$ podr\'{\i}an pertenecer a distintos elementos de $\mathscr C$, pero al ser $\mathscr C$ linealmente ordenado, alguno de estos dos elementos contiene al otro y por lo tanto contiene a ambos $u,v$), y por lo tanto existe un $(u,v)$-camino en $H$, camino que tambi\'en est\'a presente en $\tilde{H}$ y por lo tanto $\tilde{H}$ es conexa. Por otra parte, cualquier ciclo en $\tilde{H}$ tambi\'en es un ciclo en alg\'un $H\in\mathscr C$ (de entrada, cada arista del ciclo podr\'{\i}a pertenecer a alg\'un elemento distinto de $\mathscr C$, pero estos son s\'olo una cantidad finita y, al ser $\mathscr C$ linealmente ordenado, necesariamente alguno de estos elementos contiene a todos los dem\'as y por lo tanto a todo el ciclo), lo cual contradir\'{\i}a que $H$ es un \'arbol. Por lo tanto $\tilde{H}$ es ac\'{\i}clica.

As\'{\i}, $\tilde{H}\in\mathbb P$ es una cota inferior para $\mathscr C$. Por lo tanto $\mathbb P$ satisface las hip\'otesis del lema de Zorn, de modo que existe un elemento minimal $T\in\mathbb P$. Entonces $T$ es una subgr\'afica de $G$ que es un \'arbol. Si $T$ no fuera una subgr\'afica generadora, entonces podr\'{\i}amos encontrar un v\'ertice $v$ de $G$ que no pertenece a $T$; sin embargo, al ser $G$ conexa, debe de haber un $(u,v)$-camino tal que $v$ es v\'ertice de $T$. Sin p\'erdida de generalidad (recortando el camino de ser necesario) podemos suponer que ning\'un otro v\'ertice recorrido por este camino pertenece a $T$, por lo que, al a\~nadir este camino a la subgr\'afica $T$ para obtener $T'$, tendr\'{\i}amos que $T'\in\mathbb P$ y $T'\leq T$, contradiciendo que $T$ es minimal.

\hfill$\Box$

El Teorema~\ref{teor:conexaimplicaarbolgen} constituye un ejemplo bastante paradigm\'atico del fen\'omeno consistente en que resultados de combinatoria finita sigan cumpli\'endose en el caso infinito, pero las demostraciones respectivas (de los casos finito e infinito) sean completamente diferentes. En particular, la sospecha (debido a la utilizaci\'on del lema de Zorn) de que el Teorema~\ref{teor:conexaimplicaarbolgen} requiere del axioma de elecci\'on no es infundada, como lo podemos ver en el Ejemplo~\ref{ex:sinarbol} del ap\'endice a este art\'{\i}culo.

En el caso finito, el Teorema~\ref{teor:conexaimplicaarbolgen} suele utilizarse en t\'andem con el Teorema~\ref{teor:arboleigualvmasuno} para concluir, como corolario, que toda gr\'afica conexa $G$ satisface $|E|+1\geq |V|$ (pues la igualdad se cumple en cualquier \'arbol generador de $G$, y $G$ tiene por lo menos tantas aristas como cualquiera de sus \'arboles generadores). En el caso infinito esto tambi\'en es cierto, pero dista mucho de ser interesante, pues, si $G=(V,E)$ es una gr\'afica infinita y conexa, entonces (al ser conexa) debe tenerse que todo v\'ertice de $G$ es de grado al menos uno, con lo cual $|V|=|E|=|E|+1$ por el Corolario~\ref{inf-v=e}.

Cerramos esta secci\'on discutiendo brevemente la cuesti\'on relativa al conteo de la cantidad de \'arboles generadores que admite una gr\'afica dada. El n\'umero $\tau(G)$ de una gr\'afica $G$ se define (de una manera que funciona a la perfecci\'on tanto en el caso finito como en el infinito) como la cardinalidad del conjunto de todos los \'arboles generadores de $G$. Recordemos que, si $e$ es una arista de $G$, entonces $G-e$ es la gr\'afica que resulta de eliminar al elemento $e$ del conjunto de aristas de $G$ (en otras palabras, el \'unico cambio respecto de $G$ es que los dos extremos de $e$ ya no son adyacentes en $G-e$); por otra parte, recu\'erdese tambi\'en que $G\cdot e$ es la gr\'afica que resulta de {\em contraer} la arista $e$ (identificando ambos extremos de $e$ para que se conviertan en un solo v\'ertice). El siguiente resultado cl\'asico en el estudio de las gr\'aficas finitas puede trasplantarse sin mayor problema al caso infinito.

\begin{thm}\label{teor:tau}
Si $G$ es una gr\'afica, y $e$ es una arista de $G$, entonces $\tau(G)=\tau(G-e)+\tau(G\cdot e)$.
\end{thm}
    
\dem 
Basta notar que el conjunto de \'arboles generadores de $G$ se puede partir en dos subconjuntos disjuntos, a saber, el conjunto de aquellos \'arboles generadores que contienen a la arista $e$, y el de aquellos que no la contienen. El segundo conjunto coincide exactamente con el conjunto de \'arboles generadores de $G-e$; por su parte, el primer conjunto se puede poner en una biyecci\'on obvia (mapeando al \'arbol $T$ hacia $T\cdot e$) con el conjunto de \'arboles generadores de $G\cdot e$.

\hfill$\Box$

A diferencia del caso finito, en el caso infinito el Teorema~\ref{teor:tau} no nos resulta demasiado \'util para calcular los n\'umeros $\tau(G)$, ya que, en general, si $G$ tiene una cantidad infinita de aristas entonces tanto $G-e$ como $G\cdot e$ tienen conjuntos de aristas de exactamente la misma cardinalidad; esto nos impide realizar el procedimiento, que s\'{\i} funciona en el caso finito, por inducci\'on sobre el n\'umero de aristas. En general, cualquier gr\'afica $G$ con una cantidad infinita $\kappa$ de aristas t\'{\i}picamente tiene suficiente variabilidad como para poder definir $2^\kappa$ distintos \'arboles generadores en ella. A continuaci\'on demostramos esta afirmaci\'on en el caso particular m\'as f\'acil de manejar, a saber, el de las gr\'aficas completas.
 
\begin{thm}
Sea $\lambda$ un cardinal infinito. Entonces, $\tau(K_\lambda)=2^\lambda$.
\end{thm}
    
\dem
Sean $V$ el conjunto de v\'ertices, y $E$ el conjunto de aristas, de la gr\'afica completa $K_\lambda$. Desde luego, $|V|=\lambda$; tambi\'en es inmediato ver que $|E|=|[V]^2|=[\lambda]^2=\lambda$. N\'otese que podemos identificar a cada \'arbol generador de una gr\'afica con un subconjunto de su conjunto de aristas (aquellas que pertenecen al \'arbol, no habiendo ambig\"uedad en el conjunto de v\'ertices, pues este debe ser, por definici\'on, el conjunto de v\'ertices de la gr\'afica original). Esto nos proporciona una funci\'on inyectiva que manda, en este caso particular, cada \'arbol generador de $K_\lambda$ en un elemento de $\wp(E)$, y por lo tanto $\tau(K_\lambda)\leq 2^{|E|}=2^\lambda$.

\begin{figure}
\begin{center}
\includegraphics[scale=.4]{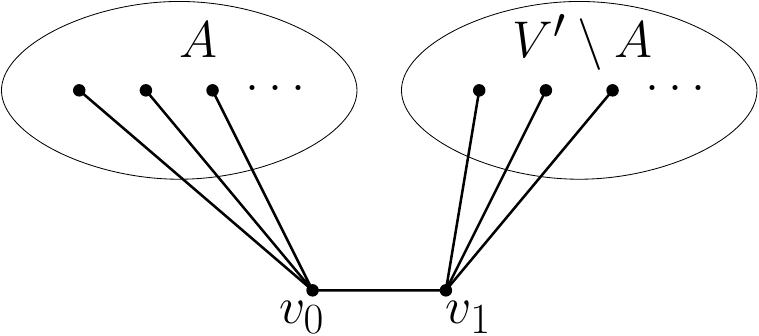}
\end{center}
\caption{El \'arbol generador $T_A$ de $K_\lambda$, para un subconjunto $A\subseteq V'=V\setminus\{v_0,v_1\}$.}
\label{fig:arbolgenerador}
\end{figure}

Procedemos ahora a demostrar la desigualdad en el sentido contrario. Fijemos dos v\'ertices $v_0,v_1\in V$, y definamos $V'=V\setminus\{v_0,v_1\}$. N\'otese que $|V'|=|V|=\lambda$. Para cada subconjunto $A\subseteq V'$, definimos la subgr\'afica $T_A$ de $K_\lambda$ como aquella subgr\'afica generadora (es decir, su conjunto de v\'ertices es $V$ mismo) cuyo conjunto de aristas viene dado por
$$
\{\{v_0,v_1\}\}\cup\{\{v_0,v\}\big|v\in A\}\cup\{\{v_1,v\}\big|v\in V'\setminus A\}.
$$
La subgr\'afica $T_A$ de $K_\lambda$ se muestra pict\'oricamente en la Figura~\ref{fig:arbolgenerador}, y no es dif\'{\i}cil demostrar formalmente que $T_A$ es, de hecho, un \'arbol generador de $K_\lambda$. M\'as a\'un, si $A,B\subseteq V'$ y $A\neq B$, entonces $T_A\neq T_B$ (pues un v\'ertice $v\in A\setminus B$ indica que la arista $\{v_0,v\}$ pertenece a $T_A$ pero no a $T_B$, y sim\'etricamente para cualquier v\'ertice en $B\setminus A$). Por lo tanto, la asignaci\'on $A\longmapsto T_A$ nos proporciona una funci\'on inyectiva desde $\wp(V')$ hacia el conjunto de \'arboles generadores de $K_\lambda$. La conclusi\'on es que $\tau(K_\lambda)\geq 2^{|V'|}=2^\lambda$, y hemos terminado.

\hfill$\Box$

\section{Coloraciones y n\'umero crom\'atico}

En esta secci\'on analizaremos algunos fen\'omenos que surgen de con\-si\-de\-rar particiones del conjunto de v\'ertices de una gr\'afica.

\begin{defi}
Dada una gr\'afica $G=(V,E)$, diremos que una {\em $\kappa$-coloraci\'on} de los v\'ertices de $V$ es alguno de los siguientes dos objetos:
\begin{enumerate}
\item una {\em partici\'on de $V$}, es decir, una familia (indexada) de subconjuntos $\{C_\alpha\subseteq V\big|\alpha<\kappa\}$ disjuntos por pares tales que $V=\bigcup_{\alpha<\kappa} C_\alpha$, o bien
\item una funci\'on $c:V\longrightarrow\kappa$.
\end{enumerate}
Diremos que la coloraci\'on es {\em apropiada} si separa v\'ertices adyacentes: esto es, una partici\'on $\{C_\alpha\big|\alpha<\kappa\}$ es apropiada si no hay dos elementos de $C_\alpha$ que sean adyacentes (en otras palabras, cada $C_\alpha$ es un {\em conjunto independiente}); equivalentemente, una funci\'on $c:V\longrightarrow\kappa$ es apropiada si para cualesquiera dos $u,v\in V$ que son adyacentes, debe tenerse que $c(u)\neq c(v)$.
\end{defi}

Desde luego que hay una traducci\'on obvia entre las dos posibles definiciones proporcionadas arriba: la funci\'on $c:V\longrightarrow\kappa$ da pie a la partici\'on $\{C_\alpha\big|\alpha<\kappa\}$ definida por $C_\alpha=c^{-1}[\{\alpha\}]$; rec\'{\i}procamente, a la partici\'on $\{C_\alpha\big|\alpha<\kappa\}$ de $V$ se le puede asociar la funci\'on $c:V\longrightarrow\kappa$ definida haciendo que, para todo $v\in V$, $c(v)$ sea el \'unico $\alpha$ tal que $v\in C_\alpha$. N\'otese que, crucialmente, la traducci\'on entre el lenguaje de las particiones y el de las funciones preserva la noci\'on de que una coloraci\'on sea, o deje de ser, apropiada. Es por ello que, a lo largo del presente art\'{\i}culo, utilizaremos la pa\-la\-bra ``coloraci\'on'' para designar a cualquiera de estos dos objetos matem\'aticos, sin hacer distinci\'on entre ellos (en cierto sentido, estos dos conceptos, si bien no son literalmente equivalentes uno con el otro, s\'{\i} lo son en alg\'un sentido m\'as metaf\'orico).

De manera puramente intuitiva, informal y extramatem\'atica, la pa\-la\-bra ``coloraci\'on'' proviene justamente de que nos estamos imaginando asignar colores a los v\'ertices de nuestra gr\'afica: al v\'ertice $v$ le asignamos el color $c(v)$ si estamos pensando en t\'erminos de funciones; equivalentemente, el conjunto $C_\alpha$ representa a todos los v\'ertices que recibieron el color $\alpha$, cuando pensamos en t\'erminos de particiones.

De particular importancia es el llamado ``n\'umero crom\'atico'' de una gr\'afica, ampliamente estudiado (tanto en el caso finito como en el infinito), y definido a continuaci\'on.

\begin{defi}
Sea $G$ una gr\'afica.
\begin{enumerate}
\item Dado un n\'umero cardinal $\kappa$, diremos que $G$ es {\em $\kappa$-crom\'atica} si admite una $\kappa$-coloraci\'on apropiada.
\item Definimos el {\em n\'umero crom\'atico} de $G$, denotado por $\chi(G)$, como el m\'{\i}nimo $\kappa$ tal que $G$ es $\kappa$-crom\'atica.
\end{enumerate}
\end{defi}

Siempre que $\kappa\leq\kappa'$, cualquier $\kappa$-coloraci\'on puede pensarse como una $\kappa'$-coloraci\'on (haciendo $C_\alpha=\varnothing$ siempre que $\kappa<\alpha<\kappa'$, si se piensa en la coloraci\'on como una partici\'on; alternativamente, considerando el codominio de la funci\'on $c$ como el conjunto $\kappa'$ en lugar de $\kappa$, si pensamos en las coloraciones como funciones); la coloraci\'on es apropiada al pensarla como una $\kappa'$-coloraci\'on si y s\'olo si ya lo era cuando se le conceb\'{\i}a como una $\kappa$-coloraci\'on. Por lo tanto, toda gr\'afica $\kappa$-crom\'atica es tambi\'en $\kappa'$-crom\'atica siempre que $\kappa\leq\kappa'$, y el n\'umero $\chi(G)$ se\~nala precisamente el lugar donde se ubica esa frontera entre admitir o no admitir coloraciones apropiadas con cierto n\'umero de colores.

Una partici\'on en singuletes (alternativamente, una funci\'on $c$ inyectiva) testifica que toda gr\'afica $G=(V,E)$ admite una $|V|$-coloraci\'on apropiada; en particular, el n\'umero $\chi(G)$ siempre existe (al menos, cuando se supone el axioma de elecci\'on; v\'ease la Proposici\'on~\ref{chinodefinido} para un ejemplo, sin el axioma de elecci\'on, de una gr\'afica cuyo n\'umero crom\'atico no est\'a bien definido). 

El axioma de elecci\'on tambi\'en resulta relevante para el concepto de la {\it compacidad} del n\'umero crom\'atico. Para entender a cabalidad esta afirmaci\'on, recordemos que el teorema de Tychonoff, uno de tantos enunciados equivalentes al axioma de elecci\'on, establece que el producto arbitrario de espacios topol\'ogicos compactos (equipado, desde luego, con la topolog\'{\i}a producto) es \'el mismo tambi\'en compacto (una demostraci\'on --utilizando, por supuesto, el axioma de elecci\'on-- del teorema de Tychonoff puede encontrarse en~\cite[Theorem 37.3]{munkres}). Recordemos, por otra parte, que un espacio topol\'ogico es compacto si y s\'olo si toda familia de subconjuntos cerrados con la PIF (propiedad de la intersecci\'on finita, que consiste en que cada intersecci\'on de una cantidad finita de elementos de la familia es no vac\'{\i}a) tiene una intersecci\'on no vac\'{\i}a (ver, por ejemplo, \cite[Theorem 26.9]{munkres}). El siguiente teorema es uno de los m\'as importantes, tanto en el estudio de las gr\'aficas infinitas, como en el marco de la teor\'{\i}a de conjuntos sin el axioma de elecci\'on.

\begin{thm}[Erd\H{o}s--de Bruijn]\label{erdos-debruijn}
Sea $G=(V,E)$ una gr\'afica infinita, y sea $k\in\mathbb N$ tal que toda subgr\'afica finita de $G$ admite una $k$-coloraci\'on apropiada. Entonces, $G$ misma admite una $k$-coloraci\'on apropiada.
\end{thm}

\dem
Consideremos el espacio topol\'ogico $X=\prod_{v\in V}\{1,\ldots,k\}$ (equipado con la topolog\'{\i}a producto que emana de equipar cada una de las copias de $\{1,\ldots,k\}$ con la topolog\'{\i}a discreta). Dado que cada uno de los espacios $\{1,\ldots,k\}$ es compacto Hausdorff, entonces $X$ tambi\'en lo es por el teorema de Tychonoff. Note que cada elemento de $X$, al ser formalmente una  funci\'on $f:V\longrightarrow\{1,\ldots,k\}$, se puede pensar como una coloraci\'on (no necesariamente apropiada) de los v\'ertices de la gr\'afica $G$. Ahora considere, para cada arista $e=\{v,u\}\in E$,
$$
F_e=\{c\in X\big|c(v)\neq c(u)\}.
$$
No es dif\'{\i}cil ver que cada conjunto $F_e$ es cerrado. Adem\'as, la familia $\{F_e\big|e\in E\}$ tiene la PIF, pues si $E'\subseteq E$ es un subconjunto finito, y $Y=\bigcup_{e\in E'}e=\{v\in V\big|v\mathrm{\ es\ extremo\ de\ alg\acute{u}n\ }e\in E'\}$, entonces
$$
\bigcap_{e\in E'}F_e=\left\{c\in X\big|\mathrm{La\ restricci\acute{o}n\ de\ }c\mathrm{\ a\ }Y\mathrm{\ es\ apropiada}\right\}.
$$
Este \'ultimo conjunto es no vac\'{\i}o pues, por hip\'otesis, cada subgr\'afica finita (en este caso, en particular la subgr\'afica inducida por $Y$) admite una coloraci\'on apropiada (y extendiendo esta coloraci\'on ar\-bi\-tra\-ria\-men\-te, sin importar si la coloraci\'on es apropiada o no fuera de $Y$, ob\-te\-ne\-mos un elemento $c\in\bigcap_{e\in E'}F_e$). De esta forma, por la compacidad de $X$, existe un elemento $c\in\bigcap_{e\in E}F_e$; es decir, un $c:V\longrightarrow\{1,\ldots,k\}$ tal que para cada arista $e=\{u,v\}\in E$, satisface $c(u)\neq c(v)$. En otras palabras, $c$ es una $k$-coloraci\'on apropiada para $G$, y hemos terminado.

\hfill$\Box$

Como consecuencia del Teorema~\ref{erdos-debruijn}, si $\chi(H)\leq k$ para toda sub\-gr\'a\-fi\-ca finita $H$ de $G$, entonces necesariamente debe tenerse que $\chi(G)\leq k$. En el pasado, se realizaron esfuerzos intensos por determinar exactamente cu\'anto del axioma de elecci\'on se utiliza en la demostraci\'on anterior. En dicha demostraci\'on, utilizamos el caso particular del teorema de Tychonoff en el cual se consideran productos topol\'ogicos de espacios Hausdorff. Ahora bien, como lo hemos dicho arriba, la versi\'on completamente general del teorema de Tychonoff es equivalente al axioma de elecci\'on; mientras tanto, la restricci\'on de este teorema que versa \'unicamente sobre espacios Hausdorff es estrictamente m\'as d\'ebil; y es, de hecho, equivalente a una versi\'on d\'ebil del axioma de elecci\'on que ha sido abundantemente estudiada, que lleva por nombre {\it teorema del Ideal Primo Booleano}. Es de hecho posible demostrar que el teorema de Erd\H{o}s--de Bruijn es equivalente al teorema del Ideal Primo Booleano. En un par de ocasiones m\'as adelante tendremos oportunidad de utilizar nuevamente este teorema, en su encarnaci\'on como teorema de Tychonoff para espacios topol\'ogicos Hausdorff.

\section{Emparejamientos en gr\'aficas}

En esta secci\'on exploramos los emparejamientos, as\'{\i} como varias relaciones entre estos y el material tratado en las secciones anteriores.

\begin{defi}
Sea $G=(E,V)$ una gr\'afica.
\begin{enumerate}
\item Un {\em emparejamiento} de $G$ es un conjunto $M\subseteq E$ tal que todo v\'ertice $v\in V$ incide en a lo m\'as una arista contenida en $M$ (e\-qui\-va\-len\-te\-men\-te, si pensamos en los elementos de $M$ como conjuntos --al ser aristas, son \'estas conjuntos de dos v\'ertices--, entonces que $M$ sea un emparejamiento es equivalente a decir que $M$ sea una familia de conjuntos disjuntos dos a dos).
\item Un emparejamiento $M$ es {\em maximal} si no existe otro emparejamiento $M'$ tal que $M\subsetneq M'$.
\item Un emparejamiento $M$ es {\em m\'aximo} si no existe otro emparejamiento $M'$ tal que $|M|<|M'|$.
\item Un emparejamiento es {\em perfecto} si todo v\'ertice $v\in V$ incide en alg\'un elemento de $M$ (pensando en los elementos de $M$ como conjuntos, equivalentemente podemos decir que $M$ es un emparejamiento perfecto si y s\'olo si $V=\bigcup_{e\in M}e$, es decir, $M$ es una partici\'on de $V$ en subconjuntos de dos elementos).
\end{enumerate}
\end{defi}

Dado un emparejamiento $M$, a los v\'ertices en los cuales incide alg\'un elemento de $M$ se les llama {\em $M$-emparejados}, o simplemente {\em emparejados} si el emparejamiento $M$ est\'a claro por el contexto. N\'otese que toda gr\'afica cuenta con al menos un emparejamiento (el conjunto vac\'{\i}o, en el peor de los casos); de ah\'{\i} tambi\'en que toda gr\'afica cuente con un emparejamiento maximal: en el caso finito, basta ir a\~nadiendo aristas una por una a un emparejamiento dado, hasta que ya no se pueda m\'as, para obtener un emparejamiento maximal; en el caso infinito basta llevar a cabo una aplicaci\'on est\'andar del lema de Zorn. Tambi\'en es cierto que toda gr\'afica cuenta con un emparejamiento m\'aximo; en el caso finito, esto puede verse con facilidad tomando, de entre todos los emparejamientos posibles --tan s\'olo hay una cantidad finita de ellos--, alguno con la m\'axima cardinalidad posible de elementos. En el caso infinito no es tan inmediato demostrar la existencia de emparejamientos m\'aximos; sin embargo, es el caso que estos existen, aunque en este art\'{\i}culo omitimos la demostraci\'on\footnote{La demostraci\'on corresponde a un trabajo, a\'un sin publicar, del primer autor en colaboraci\'on con Enrique Reyes.}.

Por otra parte, no es el caso que toda gr\'afica admita un emparejamiento perfecto, ni siquiera en el caso finito: cualquier gr\'afica con una cantidad impar de v\'ertices constituye un contraejemplo. Dado que un emparejamiento $M$ empareja a $2|M|$ v\'ertices, y el m\'aximo posible de v\'ertices a emparejar es $|V|$, no es dif\'{\i}cil ver que todo emparejamiento perfecto debe ser m\'aximo; tampoco es dif\'{\i}cil ver que todo emparejamiento perfecto es maximal. En el caso finito, un emparejamiento m\'aximo necesariamente debe de ser maximal; no as\'{\i} en el caso infinito (pi\'ensese en un emparejamiento m\'aximo infinito $M$ y t\'omese $e\in M$; entonces $M\setminus\{e\}$ tiene la misma cardinalidad de $M$ y, por lo tanto, sigue siendo m\'aximo, pero no es maximal).

A continuaci\'on trabajamos para indagar hasta qu\'e punto se puede generalizar un teorema conocido como {\em de Berge}. Para ello introducimos la siguiente definici\'on.

\begin{defi}
Sea $G=(V,E)$ una gráfica, y sea $M\subseteq E$ un emparejamiento de $G$.
\begin{enumerate}
\item Un {\em camino bi-infinito} es una sucesi\'on de v\'ertices, indexada por $\mathbb Z$, $(u_n\big|n\in\mathbb Z)$, tal que $u_n$ es adyacente a $u_{n+1}$ para todo $n\in\mathbb Z$.
\item Un {\em camino infinito en una direcci\'on} es una sucesi\'on de v\'ertices, indexada por $\mathbb N\cup\{0\}$, $(u_0,u_1,\ldots,u_n,\ldots)$, tal que $u_n$ es adyacente a $u_{n+1}$ para todo $n\in\mathbb N\cup\{0\}$.
\item Un {\em camino} es o bien un camino finito (tal y como se define al inicio de la secci\'on~\ref{sect:arboles}), o bien un camino infinito en una direcci\'on, o bien un camino bi-infinito.
\item Decimos que un camino $(u_n\big|n\in I)$ (en donde $I=\{0,\ldots,n\}$ para alg\'un $n$ si el camino es finito, $I=\mathbb N\cup\{0\}$ si el camino es infinito en una direcci\'on, o bien $I=\mathbb Z$ si el camino es bi-infinito) es {\em $M$-alternante} si sus aristas se encuentran alternadamente en $M$ y en $E\setminus M$, es decir, para cada $i\in I$ se tiene que $\{u_i,u_{i+1}\}\in M$ si y s\'olo si $\{u_{i+1},u_{i+2}\}\in E\setminus M$.
\item Decimos que un camino es {\em $M$-aumentante} si, adem\'as de ser $M$-alternante, sus v\'ertices inicial y final (en caso de que existan, es decir, en caso de que el camino sea finito, o bien \'unicamente su v\'ertice inicial, en caso de que el camino sea infinito en una direcci\'on) no est\'an $M$-emparejados.
\end{enumerate}
\end{defi}

El {\em teorema de Berge}, cl\'asico en el caso de las gr\'aficas finitas, afirma que un emparejamiento $M$ en una gr\'afica $G$ es m\'aximo si y s\'olo si $G$ no admite ning\'un camino $M$-aumentante. En el caso infinito, esto ya no se cumple, como se puede apreciar con el ejemplo de la Figura~\ref{fig:cam-infinito}, en el cual el emparejamiento $M=\{\{v_{2n},v_{2n+1}\}\big|n\in\mathbb N\}$, que es m\'aximo (al ser de cardinalidad $\aleph_0$, que es tambi\'en la cardinalidad del conjunto de aristas y por lo tanto es la m\'axima cardinalidad alcanzable por cualquier emparejamiento de nuestra gr\'afica) y, sin embargo, el camino $(v_1,v_2,\ldots,v_n,\ldots)$ satisface que su v\'ertice inicial $v_1$ no est\'a $M$-emparejado, y por lo tanto dicho camino es $M$-aumentante. Una de las implicaciones del teorema, sin embargo, s\'{\i} se cumple, como lo mostramos a continuaci\'on.

\begin{thm}
Si $G$ es una gr\'afica, y $M$ es un emparejamiento en $G$ tal que $G$ no admite caminos $M$-aumentantes, entonces $M$ es un emparejamiento máximo.
\end{thm}

\dem
Supongamos que $M$ no es máximo, y tomemos un emparejamiento $M'$ que s\'{\i} sea m\'aximo. De esta forma, debe tenerse que $|M|<|M'|$. Ahora, consideremos la gr\'afica $H$ inducida por la diferencia sim\'etrica $M\bigtriangleup M'$; note que todo v\'ertice de $H$, al ser incidente a lo más con una arista de $M$ y a lo m\'as con una de $M'$, debe de tener grado igual a $1$ o $2$. Esto implica que cada componente conexa de $H$ es un camino (finito o infinito, y en caso de ser finito, es posible que sea un ciclo); m\'as a\'un, cada uno de estos caminos es alternante entre $M$ y $M'$, pues todos los v\'ertices de $H$ son o bien $M$-emparejados o bien $M'$-emparejados, pero no ambas. En particular, ninguno de estos caminos puede ser un ciclo impar.

As\'{\i}, cada componente conexa de $H$ es o bien un ciclo par, o bien un camino infinito (ya sea infinito en una direcci\'on, o bi-infinito), o bien un camino finito que no es un ciclo. En los dos primeros casos mencionados, la componente conexa en cuesti\'on tiene tantas aristas pertenecientes a $M$ como aristas pertenecientes a $M'$. Como $|M|<|M'|$, debe de haber por lo menos una (de hecho, una infinidad de) componente conexa que tiene m\'as aristas de $M'$ que de $M$. La \'unica opci\'on para tal componente es que sea un camino finito, que no es un ciclo, de longitud impar, y cuyos v\'ertices inicial y final est\'en $M'$-emparejados pero no $M$-emparejados. En otras palabras, esta componente conexa es un camino $M$-aumentante.

\hfill$\Box$

De particular importancia en el estudio de la teor\'{\i}a de gr\'aficas son las llamadas {\em gr\'aficas bipartitas}, que son aquellas cuyo n\'umero crom\'atico es igual a $2$. Pict\'oricamente, este tipo de gr\'aficas se suelen dibujar fijando una $2$-coloraci\'on apropiada y colocando todos los v\'ertices de un color en la parte de abajo, y todos los v\'ertices del otro color en la parte de arriba, de tal forma que todas las aristas conecten v\'ertices en la parte superior del dibujo con v\'ertices en la parte inferior del mismo. En estos casos, normalmente a la coloraci\'on de la gr\'afica en dos colores, cuando se le piensa como una partici\'on del conjunto de v\'ertices, se dice que es una {\em bipartici\'on de la gr\'afica}. A continuaci\'on discutimos un poco acerca de gr\'aficas bipartitas infinitas y los criterios que nos puedan garantizar la existencia de ciertos emparejamientos, en especial perfectos, en dichas gr\'aficas. El primer resultado a lo largo de estas l\'{\i}neas es una versi\'on concreta de un resultado abstracto acerca de cardinalidades.

\begin{thm}[Cantor--Bernstein]\label{cantor-bernstein}
Sea $G$ una gr\'afica bipartita con bipartici\'on $(X,Y)$. Si existen dos emparejamientos, uno de los cuales empareja a todos los v\'ertices de $X$ y el otro que empareje a todos los v\'ertices de $Y$, entonces $G$ admite un emparejamiento perfecto.
\end{thm}

\dem
Sea $M$ un emparejamiento que empareja a todos los v\'ertices de $X$, y sea $M'$ un emparejamiento que empareja a todos los v\'ertices de $Y$. Consideremos la subgr\'afica $H=(X\cup Y,M\cup M')$ de $G$ (que es justamente la subgr\'afica que resulta de ignorar todas las aristas de $G$, excepto aquellas que pertenecen a $M$ y a $M'$). Al ser $M$ y $M'$ emparejamientos, todo v\'ertice en $H$ tiene grado $1$ o $2$; como consecuencia de ello, cada componente conexa de $H$ es un camino: un camino finito o un ciclo finito o un camino bi-infinito, o bien un camino infinito en una direcci\'on. Construiremos el emparejamiento $M''$ como sigue: para cada componente conexa de $H$ que sea infinita (ya sea que sea un camino infinito en una direcci\'on, o un camino bi-infinito), tomamos alternadamente a las aristas de dicho camino (en el caso de un camino infinito en una direcci\'on, lo hacemos comenzando por la primer arista del camino) y coloc\'andolas en $M''$. Ahora, para cada componente conexa finita, es posible argumentar que debe de tratarse de un ciclo de longitud par, por lo que tambi\'en es posible tomar alternadamente la mitad de las aristas del ciclo y colocarlas en el conjunto $M''$. Se verifica entonces que $M''$ es un emparejamiento (en $H$, y por lo tanto tambi\'en en $G$) que empareja a todos los v\'ertices de $G$.

\hfill$\Box$

Un c\'elebre criterio, en el caso finito, para la existencia de emparejamientos, es conocido por diversos autores como el {\it teorema del ma\-tri\-mo\-nio de Hall}~\cite{hall-marriage} y establece que, en una gr\'afica bipartita finita con bipartici\'on $(X,Y)$, existe un emparejamiento que empareja a todos los v\'ertices de $X$ si y s\'olo si para todo $S\subseteq X$, se cumple que $|N(S)|\geq |S|$, en donde $N(S)$ denota al conjunto de todos los v\'ertices que son adyacentes a alg\'un elemento de $S$ (es decir, $N(S)$ es el conjunto de vecinos de $S$). La raz\'on del nombre ``teorema del matrimonio'' hace alusi\'on a que, en cierta sociedad heteronormativa en la cual hipot\'eticamente quisi\'eramos generar matrimonios entre los hombres y las mujeres, re\-pre\-sen\-ta\-dos respectivamente por $X$ y $Y$, en donde una arista entre un elemento de $X$ y uno de $Y$ indica que ese hombre y esa mujer con\-si\-de\-ra\-r\'{\i}\-an aceptable casarse, el teorema nos proporciona un criterio para determinar si es posible casar ya sea a todos los hombres, o a todas las mujeres de esa sociedad. Si bien el nombre del teorema es, hoy por hoy, a todas luces obsoleto, no lo es as\'{\i} el estudio matem\'atico del mismo teorema. El an\'alogo infinito de este teorema, como vemos a con\-ti\-nua\-ci\'on, ya no es una equivalencia, aunque lo sigue siendo para una clase restringida de gr\'aficas infinitas; lo que motiva las siguientes definici\'on y teorema.

\begin{defi}
Una gr\'afica $G$ es {\em localmente finita} si el grado de cada uno de sus v\'ertices es finito.
\end{defi}

Es inmediato notar que, en particular, toda gr\'afica finita es localmente finita.

\begin{thm}\label{teor:hall-marriage}
Sea $G$ una gr\'afica bipartita con bipartición $(X,Y)$. Entonces:
\begin{enumerate}
\item Si $G$ admite un emparejamiento que empareja a cada elemento de $X$, entonces todo $S\subseteq X$ satisface $|N(S)|\geq |S|$.
\item Si adem\'as $G$ es localmente finita, y todo $S\subseteq X$ satisface que $|N(S)|\geq|S|$, entonces $G$ admite un emparejamiento que empareja a todos los elementos de $X$.
\end{enumerate}
\end{thm}

\dem\hfill
\begin{enumerate}
\item Sea $M$ un emparejamiento que empareja a todos los elementos de $X$, y sea $S\subseteq X$. Entonces, cada $v\in S$ es incidente en alguna arista $e_v\in M$; definamos $f(v)$ como el extremo de $e_v$ que es distinto de $v$, para todo $v\in S$, y notemos que $f(v)$ siempre es un elemento del conjunto $N(S)$. Entonces, la funci\'on $f:S\longrightarrow N(S)$ es inyectiva, lo cual muestra que $|S|\leq |N(S)|$.
\item Supongamos ahora que $G=(V,E)$ es localmente finita y que satisface la hip\'otesis de que $|S|\leq|N(S)|$ para todo $S\subseteq X$. Utilizaremos una vez m\'as el teorema de Tychonoff; para ello, con\-si\-de\-re\-mos el espacio topol\'ogico $\prod_{e\in E}\{0,1\}$, el cual es compacto al serlo cada uno de sus factores (cada factor $\{0,1\}$ viene equipado con la topolog\'{\i}a discreta). Crucialmente, n\'otese que cada punto $x$ de nuestro espacio topol\'ogico es la funci\'on caracter\'{\i}stica de alg\'un subconjunto $M_x=\{e\in E\big|x(e)=1\}$ de $E$. Para cada $S\subseteq X$, sea
\begin{eqnarray*}
F_S=\bigg\{x\in\prod_{e\in E}\{0,1\} & \bigg| & M_x\mathrm{\ es\ un\ emparejamiento} \\
& & \mathrm{que\ empareja\ a\ cada\ elemento\ de\ }S\bigg\}
\end{eqnarray*}
A continuaci\'on verificaremos que cada uno de los subconjuntos $F_S$ es cerrado. Pues si tenemos alg\'un $S\subseteq X$ y alg\'un $x\notin F_S$, entonces o bien $M_x$ no es un emparejamiento (es decir, existen dos aristas $e_1,e_2\in M_x$ que tienen alguno de sus extremos en com\'un), o bien s\'{\i} lo es, pero hay alg\'un $v\in S$ que no es emparejado por $M_x$ (es decir, si $e_1,\ldots,e_n$ son las aristas incidentes en $v$, entonces cada una de estas $e_i\notin M_x$; note que la cantidad de aristas incidentes en $v$ es finita por ser $G$ localmente finita). En el primer caso hacemos $U=\{x\in\prod_{e\in E}\{0,1\}\big|x(e_1)=1=x(e_2)\}$; en el segundo caso de\-fi\-ni\-mos $U=\{x\in\prod_{e\in E}\{0,1\}\big|x(e_i)=0\mathrm{\ para\ cada\ }i\in\{1,\ldots,n\}\}$; en ambos casos $U$ es un subconjunto abierto de $\prod_{e\in E}\{0,1\}$, completamente disjunto con $F_S$, tal que $x\in U$, lo cual demuestra que el complemento del conjunto $F_S$ es abierto. Adem\'as, siempre que $S\subseteq X$ sea un conjunto finito, tenemos que $F_S$ es no vac\'{\i}o, ya que por el teorema de Hall finito siempre podemos elegir un emparejamiento en la subgr\'afica (finita) de $G$ inducida por $S\cup N(S)$ que empareje a todo elemento de $S$; la funci\'on caracter\'{\i}stica de tal emparejamiento es un elemento de $F_S$. Finalmente, note que, para $S_1,\ldots,S_k\subseteq X$, se cumple que $F_{S_1}\cap\cdots\cap F_{S_k}=F_{S_1\cup\cdots\cup S_k}$ y por lo tanto, en caso de que cada uno de los $S_i$ sea  finito, esta intersecci\'on ser\'a no vac\'{\i}a. As\'{\i}, la familia de conjuntos cerrados $\{F_S\big|S\subseteq X\mathrm{\ es\ finito}\}$ tiene la PIF, lo que implica que su intersecci\'on es no vac\'{\i}a. No es dif\'{\i}cil verificar que, si $\displaystyle{x\in\bigcap_{S\subseteq X\mathrm{\ finito}}F_S}$, entonces en particular $x\in F_{\{v\}}$ para cada $v\in X$, por lo que $M_x$ es un emparejamiento de $G$ que empareja a todo elemento de $X$.
\end{enumerate}

\hfill$\Box$

\begin{figure}
\begin{center}
\includegraphics[scale=.6]{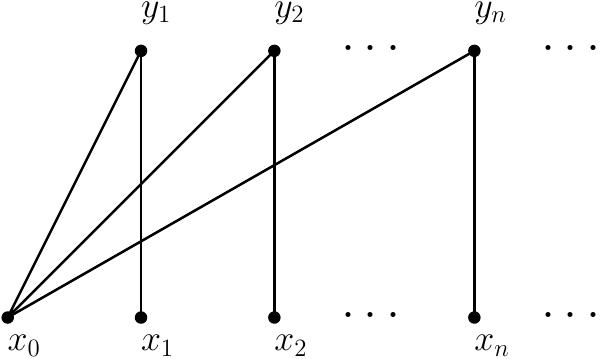}
\end{center}
\caption{Una gr\'afica bipartita que satisface el criterio del teorema de Hall pero no admite emparejamientos que emparejen a todos los v\'ertices $x_n$. Al v\'ertice $x_0$, que es susceptible de emparejarse con cualquiera de los $y_n$, se le conoce como {\it el playboy}.}
\label{fig:playboy}
\end{figure}

En el caso de las gr\'aficas que no son localmente finitas, el rec\'{\i}proco de la implicaci\'on contenida en la primera parte del Teorema~\ref{teor:hall-marriage} no se cumple, como puede verse considerando a la gr\'afica $G=(X\cup Y,E)$, en donde $X=\{x_0\}\cup\{x_n\big|n\in\mathbb N\}$ y $Y=\{y_n\big|n\in\mathbb N\}$ son disjuntos, y 
$$
E=\left\{ \{ x_0 ,y_n\}\big|n\in\mathbb{N}\right\} \cup \left\{ \{x_n,y_n\}\big|n\in\mathbb{N}\right\}.
$$
Esta gr\'afica est\'a ilustrada en la Figura~\ref{fig:playboy}. No es dif\'{\i}cil ver que, para cada $S\subseteq X$, $|N(S)|\geq|S|$ (pues $|N(S)|$ es o bien $|S|$, si $x_0\notin S$, o bien $\aleph_0$, en caso contrario); por otra parte, $G$ no admite ning\'un emparejamiento que empareje a todo elemento de $X$, pues tan pronto como $x_0$ est\'a emparejado, digamos que por medio de la arista $\{x_0,y_n\}$, de inmediato se vuelve imposible emparejar al v\'ertice $x_n$.

Si $\lambda$ es un n\'umero cardinal, una gr\'afica se dice que es {\em $\lambda$-regular} si todos sus v\'ertices tienen grado $\lambda$; una gr\'afica es {\em regular} si es $\lambda$-regular para alg\'un $\lambda$. En el contexto de las gr\'aficas finitas, una consecuencia in\-me\-dia\-ta del teorema del matrimonio de Hall\footnote{De hecho, algunos autores, por Ejemplo~\cite{bondy}, le llaman {\em teorema del matrimonio de Hall} no al enunciado del Teorema~\ref{teor:hall-marriage}, sino al enunciado que presentamos a continuaci\'on.} es que toda gr\'afica bipartita regular con bipartici\'on $(X,Y)$ debe satisfacer $|X|=|Y|$. En el caso infinito, probar que toda gr\'afica bipartita regular con bipartici\'on $(X,Y)$ satisface $|X|=|Y|$ es extremadamente sencillo, sin necesidad de utilizar ning\'un otro resultado como lema previo (y ser\'a un ejercicio interesante para el lector desarrollar tal demostraci\'on); sin embargo, como hemos resaltado anteriormente, la simple comparaci\'on de cardinalidades por lo general no nos proporciona enunciados lo suficientemente informativos al momento de estudiar gr\'aficas infinitas. El resultado que presentamos a continuaci\'on s\'{\i} es lo suficientemente informativo, pues nos proporcionar\'a la biyecci\'on entre las dos clases de la bipartici\'on no s\'olo en abstracto, sino en t\'erminos de las aristas de la gr\'afica en cuesti\'on. 

\begin{thm}
Sea $G$ una gr\'afica bipartita regular. Entonces, $G$ admite un emparejamiento perfecto.
\end{thm}

\dem
Si $G$ es localmente finita, el resultado se sigue de aplicar el teorema del matrimonio de Hall (Teorema~\ref{teor:hall-marriage}) junto con el teorema de Cantor--Bernstein (Teorema~\ref{cantor-bernstein}). De modo que podemos suponer sin p\'erdida de generalidad que $G$ es $\kappa$-regular, para alg\'un cardinal infinito $\kappa$. En este caso, note que cada componente conexa de $G$ debe de tener cardinalidad $\kappa$: dado un v\'ertice $v_0$ en dicha componente conexa, la componente es simplemente $\bigcup_{n\in\mathbb N\cup\{0\}}\{v\in V\big|d(v_0,v)=n\}$, y cada uno de los conjuntos $\{v\in V\big|d(v_0,v)=n\}$ tiene a lo m\'as $\kappa$ elementos, como es f\'acil demostrar de manera inductiva y utilizando el hecho de que $G$ es $\kappa$-regular. As\'{\i}, el emparejamiento perfecto se puede construir independientemente en cada componente conexa, lo que se traduce, en t\'erminos pr\'acticos, a que podemos suponer sin perder generalidad que $|V|=\kappa$. Bien-ordenamos $V=\{v_\alpha\big|\alpha<\kappa\}$, y recursivamente cons\-trui\-mos el emparejamiento $M=\{e_\alpha\big|\alpha<\kappa\}$, eligiendo simplemente, en el $\alpha$-\'esimo paso, y en caso de que $v_\alpha$ no sea incidente a ninguna arista en el conjunto $\{e_\xi\big|\xi<\alpha\}$, cualquier arista $e_\alpha$ incidente en $v_\alpha$ cuyo otro extremo no sea extremo de ninguna $e_\xi$ con $\xi<\alpha$ (lo cual se puede hacer pues $\{e_\xi\big|\xi<\alpha\}$ es un conjunto de cardinalidad estrictamente menor que $\kappa$ y $v_\alpha$ es un v\'ertice de grado $\kappa$). Al final del proceso, el emparejamiento $M=\{e_\alpha\big|\alpha<\kappa\}$ es perfecto. 

\hfill$\Box$

Acerc\'andonos lentamente a la culminaci\'on de esta secci\'on, mencionaremos, sin demostrar (debido a su grado de complejidad), algunos resultados de inter\'es may\'usculo desde el punto de vista hist\'orico. Para ello, habremos de enunciar la siguiente definici\'on.

\begin{defi}
Sea $G=(V,E)$ una gr\'afica.
\begin{enumerate}
\item Una {\em cubierta por v\'ertices} de $G$ es un conjunto $C\subseteq V$ tal que toda arista de $G$ incide en alg\'un elemento de $C$.
\item La cubierta por v\'ertices $C$ de $G$ se dice {\em minimal} si no existe otra cubierta por v\'ertices $C'$ tal que $C'\subsetneq C$.
\item La cubierta por v\'ertices $C$ de $G$ se dice {\em m\'inima} si no existe otra cubierta por v\'ertices $C'$ tal que $|C'|<|C|$.
\end{enumerate}
\end{defi}

No es dif\'{\i}cil ver (tanto en el caso finito, como en el infinito) que, si $M$ es cualquier emparejamiento y $C$ es cualquier cubierta, entonces $|M|\leq|C|$, pues cualquier funci\'on que mande a cada $e\in M$ a alg\'un $v\in C$ que sea extremo de $e$ (siempre hay por lo menos un extremo con esta caracter\'{\i}stica, al ser $C$ una cubierta) debe de ser una funci\'on inyectiva al ser $M$ un emparejamiento. En particular, esta desigualdad sigue cumpli\'endose cuando $M$ es un emparejamiento m\'aximo y $C$ es una cubierta m\'{\i}nima. El siguiente teorema es uno bastante cl\'asico en la teor\'{\i}a de gr\'aficas finitas.

\begin{thm}[K\"onig]
Si $G$ es una gr\'afica bipartita (finita), entonces la cardinalidad de cualquier emparejamiento m\'aximo de $G$ es igual a la cardinalidad de cualquier cubierta m\'{\i}nima de $G$.
\end{thm}

No es dif\'{\i}cil demostrar que, en el caso de una gr\'afica bipartita infinita, el enunciado anterior acerca de cardinalidades sigue siendo cierto; nuevamente, este resultado no resulta demasiado informativo pues la cardinalidad no nos proporciona una distinci\'on tan fina entre los diversos conjuntos infinitos. Una versi\'on m\'as expl\'{\i}cita del resultado anterior nos sugerir\'{\i}a buscar, m\'as que una simple igualdad de cardinalidades entre un emparejamiento m\'aximo $M$ y una cubierta m\'{\i}nima $C$, una correspondencia dada en t\'erminos de la gr\'afica --por ejemplo, la po\-si\-bi\-li\-dad de que a cada $e\in M$ le corresponda un \'unico $c\in C$ tal que $c$ sea extremo de $e$--. Este resultado m\'as fuerte, cuya demostraci\'on constituye todo un {\it tour de force} t\'ecnico, de complejidad bastante considerable, fue demostrado por R. Aharoni.

\begin{thm}[Aharoni~\cite{aharoni:konig}]
Sea $G$ una gr\'afica bipartita (finita o infinita). Entonces, existe un emparejamiento $M$, as\'{\i} como una cubierta por v\'ertices $C$, tales que $C$ consta de elegir exactamente un extremo de cada elemento de $M$.
\end{thm}

Lentamente nos encaminaremos a enunciar un resultado a\'un m\'as fuerte. Comencemos por recordar que, para una gr\'afica conexa no completa $G$, se define $\kappa(G)$, la {\it conectividad} de la gr\'afica, como el m\'{\i}nimo n\'umero de v\'ertices que debe uno eliminar de $G$ para que el resultado sea disconexo; an\'alogamente pero con aristas se define la {\it conectividad por aristas} $\kappa'(G)$. (Es posible definir tambi\'en estos par\'ametros en una gr\'afica no conexa, pero entonces no es tan inmediata la generalizaci\'on al caso infinito, por lo que en este caso \'unicamente estudiaremos los n\'umeros de conexidad en gr\'aficas conexas.) A un conjunto $S$ de v\'ertices tal que $G-S$ es disconexo, o bien a un conjunto $L$ de aristas tal que $G-L$ es disconexo, se les conoce como {\em conjuntos de corte} (si es necesario, se especifica si el conjunto de corte es por v\'ertices o por aristas, seg\'un sea el caso). El siguiente resultado es bien conocido en el caso finito; a continuaci\'on proporcionamos una demostraci\'on para el enunciado an\'alogo correspondiente a las gr\'aficas infinitas.

\begin{thm}
Si $G=(V,E)$ es una gráfica infinita, conexa, y no completa, entonces $\kappa(G)\leq\kappa'(G)\leq\delta(G)$.
\end{thm}

\dem
La desigualdad $\kappa'(G)\leq\delta(G)$ se establece eligiendo un v\'ertice $v\in V$ de grado $\delta(G)$, de modo y manera que, si $S$ es el conjunto de aristas incidentes en $v$, necesariamente $G-S$ ha de ser disconexa (como m\'{\i}nimo, el singulete $\{v\}$ es una componente conexa de $G-S$, por lo que esta \'ultima tiene por lo menos otra componente que contiene a los dem\'as v\'ertices); as\'{\i}, $S$ es un corte por aristas de cardinalidad $\delta(G)$.

Para establecer la desigualdad $\kappa(G)\leq\kappa'(G)$, comencemos por notar que esta desigualdad se cumple autom\'aticamente en caso de que $\kappa'(G)=|V|$ (ya que necesariamente, por definici\'on, se cumple siempre que $\kappa(G)\leq|V|$). As\'{\i} que podemos suponer sin p\'erdida de generalidad que $\kappa'(G)<|V|$. Tomemos un corte por aristas $L$ con $|L|=\kappa'(G)$, y sea $Z$ el conjunto de v\'ertices incidentes con alg\'un elemento de $L$ (t\'ecnicamente, $Z=\bigcup_{e\in L}e$). Note que $|Z|\leq 2|L|$ y, en particular, e\-xis\-ten v\'ertices de $G$ que no pertenecen a $Z$. Ahora, cada componente conexa de $G-L$ debe intersectar a $Z$, de modo que $Z$ admite una partición dada por sus intersecciones con las distintas componentes conexas de $G-L$. Uniendo miembros de esa partici\'on de ser necesario, podemos obtener dos conjuntos disjuntos no vac\'{\i}os $X,Y$ tales que $Z=X\cup Y$ y de modo que ninguna componente conexa de $G-L$ intersecta a ambos $X$ y $Y$ (en particular, todo camino de alg\'un elemento de $X$ a alg\'un elemento de $Y$ debe involucrar alguna arista de $L$; m\'as a\'un, todo tal camino involucra por lo menos a una arista que tiene un extremo en $X$ y el otro en $Y$).

Sea $v\in V$ un v\'ertice tal que $v\notin Z$, y supongamos sin p\'erdida de generalidad que la componente conexa de $v$ en $G-L$ intersecta a $X$ (y por lo tanto no intersecta a $Y$). No es dif\'icil verificar entonces que $X$ constituye un conjunto de corte por v\'ertices para $G$: por ejemplo, ning\'un elemento de $Y$ está conectado con $v$ en $G-X$ (de lo contrario, cualquier camino de alg\'un $y\in Y$ a $v$ puede utilizarse para construir un camino de $y$ a alg\'un $x\in X$, pero todo camino de tal naturaleza involucra alguna arista que tiene un extremo en $X$ y el otro en $Y$, y tales aristas ya no pertenecen a $G-X$). Por lo tanto, s\'olo resta demostrar que $|X|\leq|L|=\kappa'(G)$ para llegar a la conclusi\'on deseada. Ahora, si $|L|=\kappa'(G)$ es infinito, entonces tenemos que $|L|=|Z|=|X|+|Y|=\max\{|X|,|Y|\}$ y en particular $|X|\leq|L|$, por lo que hemos terminado. De modo que el resto de la demostraci\'on se enfocar\'a en el caso cuando $|L|=\kappa'(G)$ es finito. En este caso, comencemos por notar que toda arista de $L$ tiene un extremo en $X$ y el otro extremo en $Y$. De lo contrario, si hubiera por ejemplo $y,y'\in Y$ tales que $\{y,y'\}\in L$, entonces $L\setminus\{\{y,y'\}\}$ ya es un conjunto de corte por aristas en $G$ (pues cualquier camino conectando a $y$ con alg\'un $x\in X$ involucra al menos una arista con un extremo en $X$ y el otro en $Y$, por lo que este camino existe en $G-(L\setminus\{y,y'\})$ si y s\'olo si existe en $G-L$). Como $L$ es finito, tenemos que $\kappa'(G)\leq|L\setminus\{y,y'\}|=|L|-1=\kappa'(G)-1$, una contradicci\'on. El argumento para ver que $L$ no contiene aristas con ambos extremos en $X$ es enteramente similar. Por otra parte, note tambi\'en que {\it toda} arista con un extremo en $X$ y el otro en $Y$ debe pertenecer a $L$ (de lo contrario habría caminos de alg\'un elemento de $X$ a alg\'un elemento de $Y$ en $G-L$). De modo que cada elemento $x\in X$ es el extremo de por lo menos una arista de $L$; eligiendo de esta manera una arista de $L$ por cada $x\in X$ obtenemos una funci\'on inyectiva de $X$ en $L$, lo que implica que $|X|\leq|L|$.

\hfill$\Box$

Se dice que dos caminos son {\em internamente disjuntos} si \'unicamente tienen en com\'un sus v\'ertices iniciales y finales. En general, decimos que una familia de caminos $\mathscr F$ es {\em internamente disjunta} si cualesquiera dos de sus elementos son internamente disjuntos. Varios de los resultados, en el caso finito, m\'as b\'asicos acerca de conectividad que involucran caminos internamente disjuntos siguen siendo ciertos, y con exactamente la misma demostraci\'on que en el caso infinito. Por ejemplo, toda gr\'afica (finita o infinita) es 2-conexa (es decir, satisface $\kappa(G)\geq 2$) si y s\'olo si para cualesquiera dos v\'ertices, existen dos caminos internamente disjuntos que conectan a esos v\'ertices.

A continuaci\'on, enunciamos el {\it teorema de Menger}. Recu\'erdese que, dados dos conjuntos de v\'ertices $X,Y$ en una gr\'afica, un conjunto $S$ de v\'ertices es {\em $(X,Y)$-separador} si todo camino de $X$ a $Y$ intersecta a $S$ --equivalentemente, si en $G-S$ tenemos que ninguna componente conexa intersecta simult\'aneamente a $X$ y a $Y$. El {\em teorema de Menger} afirma que, para cualesquiera dos conjuntos disjuntos de v\'ertices $X,Y$, el tama\~no m\'{\i}nimo de un conjunto $(X,Y)$-separador es igual al tama\~no m\'aximo de una familia internamente disjunta de caminos de $X$ a $Y$. Este teorema se encuentra, por motivos hist\'oricos, estrechamente ligado al teorema de K\"onig: en el momento en que Menger public\'o su teorema, su demostraci\'on ten\'{\i}a el problema de asumir, sin demostraci\'on, el caso bipartito --este hueco en la demostraci\'on fue justamente llenado por K\"onig al momento de publicar el teorema que lleva su nombre--.

No es dif\'{\i}cil establecer como cierto el enunciado del teorema de Menger traducido palabra por palabra al caso infinito, pero, como ya viene siendo un tema recurrente, simplemente hablar de cardinalidades no es lo suficientemente informativo en este caso. Erd\H{o}s conjetur\'o que una versi\'on m\'as fuerte, y m\'as ``expl\'{\i}cita'', del teorema, se podr\'{\i}a demostrar. Finalmente, R. Aharoni y E. Berger lograron eventualmente establecer la conjetura de Erd\H{o}s, resultado tanto de complejidad como de importancia enormemente significativa, que a continuaci\'on e\-nun\-cia\-mos sin demostraci\'on.

\begin{thm}[Aharoni--Berger~\cite{aharoni-berger}]
Sea $G=(V,E)$ una gr\'afica infinita, y sean $X,Y\subseteq V$. Entonces, existe una familia internamente disjunta $\mathscr F$ de caminos que van de $X$ a $Y$, y existe un conjunto $S$ que es $(X,Y)$-separador, tal que $S$ consta de elegir exactamente un v\'ertice interno de cada uno de los elementos de $\mathscr F$.
\end{thm}        

Para concluir la secci\'on, mencionaremos un criterio para la existencia de emparejamientos perfectos en una gr\'afica. En el caso de una gr\'afica finita $G$, se define $o(G)$ como la cantidad de componentes conexas {\em impares} (es decir, con una cantidad impar de v\'ertices); un resultado cl\'asico y b\'asico en este contexto es que una gr\'afica finita $G$ admite un emparejamiento perfecto si y s\'olo si $o(G-S)\leq|S|$ para todo $S\subseteq V$.

En el caso de una gr\'afica infinita, la primera reacci\'on es extender las definiciones declarando que una componente infinita de $G$ cuente como una componente par (finalmente, todo conjunto de cardinalidad infinita se puede partir en dos subconjuntos de cardinalidades iguales). Sin embargo, con esta definici\'on \'unicamente se tiene una de las implicaciones con las que se cuenta en el caso finito: si $G=(V,E)$ es una gr\'afica que admite un emparejamiento perfecto, entonces para cada $S\subseteq V$ se cumple que $o(G-S)\leq|S|$. Pues para cada elemento de $o(G-S)$ se tiene por lo menos un v\'ertice que no est\'a emparejado con otro v\'ertice en la misma componente de $G-S$, por lo cual la \'unica opci\'on es que est\'e emparejado con un elemento de $S$, y esto define una funci\'on inyectiva de $o(G-S)$ en $S$. Se queda como ejercicio para el lector verificar que la gr\'afica de la Figura~\ref{fig:playboy}, misma que ya vimos que no admite un emparejamiento perfecto, es un contraejemplo para la implicaci\'on rec\'{\i}proca.

Una vez m\'as, nos topamos con un ejemplo que muestra que, en ocasiones, la clave para generalizar un resultado al caso infinito reside en la elecci\'on adecuada de las definiciones. Se dice que una gr\'afica es {\em factor-cr\'{\i}tica} si, para todo v\'ertice $v$ de $G$, la gr\'afica $G-v$ admite un emparejamiento perfecto. Note que una gr\'afica factor-cr\'{\i}tica finita debe de tener una cantidad impar de v\'ertices; en el caso infinito, esta es la definici\'on an\'aloga ``correcta'' de componente impar. Entonces, se cuenta con el siguiente resultado.

\begin{thm}[Aharoni~\cite{aharoni-matching}]
Si $G$ es cualquier gr\'afica, entonces $G$ admite un emparejamiento perfecto si y s\'olo si para todo $S\subseteq V(G)$, es posible emparejar de manera inyectiva a cada componente factor-cr\'{\i}tica de $G-S$ con los elementos de $S$.
\end{thm}

\section*{Ap\'endice}

En este ap\'endice presentamos un par de ejemplos, construidos (consistentemente) dentro de la teor\'{\i}a $\zf$ (es decir, dentro de los axiomas de Zermelo--Fraenkel, sin el axioma de elecci\'on), ejemplificando la manera en que varios de los teoremas demostrados anteriormente utilizan de manera esencial dicho axioma. Dado que no pretendemos introducir al lector a las vicisitudes t\'ecnicas de las pruebas de consistencia e independencia en $\zf$ (toda una teor\'{\i}a de complejidad avanzada, que f\'acilmente puede requerir de varios cursos, a lo largo de varios semestres, para poder devenir familiarizado con ella), nos apoyaremos de manera sustancial en la confianza que nos tenga el lector al momento de afirmar que ciertos objetos (en concreto, un conjunto de Russell) existen de manera consistente con $\zf$ --es decir, que no se puede demostrar la no existencia de dicho objeto \'unicamente con los axiomas de Zermelo--Fraenkel--. Sea como fuere, nuestras afirmaciones de consistencia vienen acompa\~nadas de la correspondiente referencia en la cual el lector interesado, y versado en las t\'ecnicas relevantes (tanto la t\'ecnica de los modelos de permutaciones de Fraenkel--Mostowski, como en las t\'ecnicas del forzamiento y los modelos sim\'etricos), pueda verificar la veracidad de nuestras afirmaciones.

Para construir nuestros ejemplos resulta crucial el concepto de un {\em conjunto de Russell}. Un conjunto de Russell es un conjunto $X$ que admite una partici\'on, $X=\bigcup_{n\in\mathbb N}P_n$, tal que cada $P_n$ es un conjunto de cardinalidad exactamente dos (un par); y de tal suerte que ning\'un subconjunto infinito de la familia $\{P_n\big|n\in\mathbb N\}$ admite una funci\'on de elecci\'on. Durante el resto de esta secci\'on, mantendremos la notaci\'on de $X$ para un conjunto de Russell, y $P_n$ para los pares que componen a dicho conjunto. Los conjuntos de Russell constituyen la contraparte formal del ejemplo de ``una infinidad de parejas de calcetines'' de B. Russell, y es posible demostrar que, consistentemente con los axiomas de $\zf$, tales conjuntos existen: por ejemplo, el conjunto de \'atomos en el segundo modelo de Fraenkel es un conjunto de Russell, v\'ease~\cite[Sec. 4.4]{jech-choice}.

\begin{ejem}\label{ex:secondfraenkel}
{\em Una gr\'afica $G$, tal que $\sum_{v\in V(G)}d_G(v)\neq 2|E(G)|$.} Comenzamos por tomar un conjunto de Russell $X$, y nuestra gr\'afica ser\'a $G=(X,E)$ en donde $E=\{P_n\big|n\in\mathbb N\}$. Note que el conjunto de aristas de $G$ es numerable (est\'a indexado por el conjunto $\mathbb N$), por lo cual $2|E|=2\aleph_0=\aleph_0$. Note tambi\'en que cada v\'ertice de $G$ tiene grado $1$, por lo que tenemos que $\sum_{v\in V}d(v)=|X|$ (m\'as formalmente, siguiendo la notaci\'on de la demostraci\'on del Teorema~\ref{teor:2e}, para cada $x\in X$ podemos definir $X_x=\{(x,y)\}$, en donde $\{x,y\}=P_n$ para alg\'un $n$; no es dif\'{\i}cil verificar que la funci\'on $x\longmapsto(x,y)$, en donde $y$ y $n$ son los \'unicos tales que $\{x,y\}=P_n$, constituye una biyecci\'on entre $X$ y $\bigcup_{x\in X}X_x$). Sin embargo, $X$ no puede ser equipotente a $\aleph_0$ --ya que, de hecho, la existencia de una biyecci\'on entre $X$ y $\mathbb N$, digamos que representada por $X=\{x_n\big|n\in\mathbb N\}$, implicar\'{\i}a que la funci\'on $n\longmapsto x_k$, en donde $k=\min\{m\in\mathbb N\big|x_m\in P_n\}$, constituya una funci\'on de elecci\'on para la familia $\{P_n\big|n\in\mathbb N\}$, una contradicci\'on.

Es incluso posible verificar que $|X|$ no s\'olo no es igual a $\aleph_0$, sino que estas dos cardinalidades ni siquiera son comparables\footnote{Cabe destacar que el enunciado ``las cardinalidades est\'an totalmente ordenadas'' es equivalente al axioma de elecci\'on.}; de esta forma, en nuestro ejemplo hemos obtenido que los cardinales $\sum_{v\in V(G)}d_G(v)$ y $2|E(G)|$ no son comparables.
\end{ejem}

\begin{figure}
\begin{center}
\includegraphics[scale=.6]{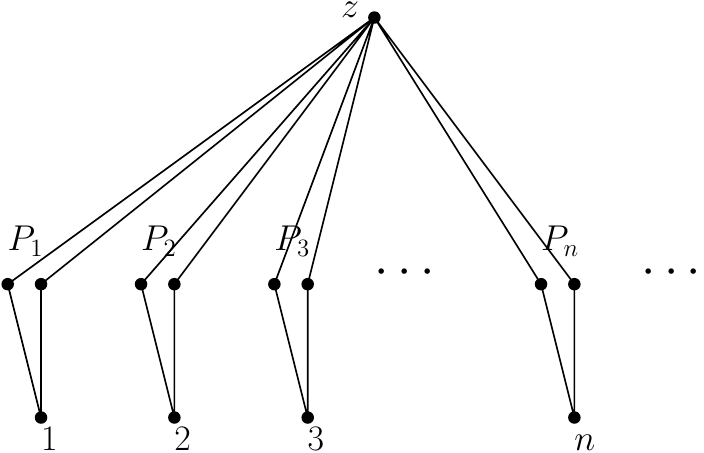}
\end{center}
\caption{Una gr\'afica conexa sin \'arboles generadores (en el mundo sin el axioma de elecci\'on).}
\label{fig:no-arbol}
\end{figure}

\begin{ejem}\label{ex:sinarbol}
{\em Una gr\'afica conexa que no admite un \'arbol generador.} Nuevamente, tomamos a nuestro conjunto de Russell $X$, tomemos un objeto $z\notin X$, y definamos la gr\'afica $G=(V,E)$ mediante
$$
V=X\cup\{z\}\cup\mathbb N,
$$
y
$$
E=\{\{z,x\}\big|x\in X\}\cup\{\{n,x\}\big|n\in\mathbb N\mathrm{\ y\ }x\in P_n\}.
$$
La gr\'afica $G$ se ilustra en la Figura~\ref{fig:no-arbol}. Es f\'acil ver que $G$ es conexa; tampoco es complicado verificar que, para cada $n\in\mathbb N$, existen exactamente dos $(n,z)$-caminos en $G$ (ambos de longitud dos, en donde los v\'ertices intermedios son respectivamente cada uno de los dos elementos de $P_n$). Por lo tanto, si $G$ admitiera un \'arbol generador $T$, la elecci\'on de un \'unico $(n,z)$-camino que pertenece a $T$ nos proporcionar\'{\i}a una manera de elegir a un \'unico elemento de $P_n$ (a saber, el \'unico que pertenece al \'unico $(n,z)$-camino en $T$). En otras palabras, a partir de un \'arbol generador de $G$ podr\'{\i}amos construir una funci\'on de elecci\'on para la familia $\{P_n\big|n\in\mathbb N\}$, lo cual es imposible; por lo tanto, $G$ carece de \'arboles generadores.
\end{ejem}

En este momento pasamos a considerar el n\'umero crom\'atico de una gr\'afica. En este contexto, el Ejemplo~\ref{ex:secondfraenkel} ilustra un fen\'omeno bastante interesante.

\begin{prop}\label{chinodefinido}
La gr\'afica del Ejemplo~\ref{ex:secondfraenkel} no tiene un n\'umero crom\'atico bien definido.
\end{prop}

\dem
En el contexto de la teor\'{\i}a de conjuntos sin el axioma de elecci\'on, supondremos que a cada conjunto $X$ se le puede asignar de manera \'unica un objeto, denotado por $|X|$, de tal forma que para cualesquiera dos conjuntos $X$ y $Y$ se satisfaga que $|X|=|Y|$ si y s\'olo si $X$ es equipotente a $Y$\footnote{En principio, nos gustar\'{\i}a que, para cada $X$, $|X|$ fuera el conjunto de todos los conjuntos equipotentes a $X$ (es decir, la clase de equivalencia de $X$, m\'odulo equipotencia). Esto no es posible, ya que dicha colecci\'on es de hecho una clase propia; sin embargo, la t\'ecnica conocida como {\it truco de Scott} nos permite definir al objeto $|X|$ de manera que tenga las propiedades deseadas. Dado que los pormenores del truco de Scott constituyen una porci\'on moderadamente avanzada de la teor\'{\i}a de conjuntos, en este art\'{\i}culo los omitimos, y simplemente suponemos que la asignaci\'on $X\longmapsto|X|$ est\'a bien definida.}.

Entonces, con la misma notaci\'on que en el Ejemplo~\ref{ex:secondfraenkel}, su\-pon\-ga\-mos que $A$ es un conjunto (de colores) y que $c:X\longrightarrow A$ es una coloraci\'on adecuada. Note que, para cada $a\in A$, el conjunto $c^{­-1}[\{a\}]$ debe de ser finito. De lo contrario, dado que a lo m\'as un elemento de cada $P_n$ puede tener color $a$ (es decir, $|c^{-1}[\{a\}]\cap P_n|\leq 1$), entonces esto nos proporcionar\'{\i}a una funci\'on de elecci\'on para alg\'un subconjunto infinito de $\{P_n\big|n\in\mathbb N\}$, lo cual es una contradicci\'on. De hecho, podemos definir la funci\'on $g:A\longrightarrow X$ dada por $g(a)=x$, en donde $x$ es tal que $c(x)=a$ y $x\in P_n$ para $n\in\mathbb N$ m\'{\i}nimo tal que $c^{­-1}[\{a\}]\cap P_n\neq\varnothing$. De esta forma, la funci\'on $g$ es inyectiva, y por lo tanto su imagen $\{g(a)\big|a\in A\}$ debe, para todo $n$ lo suficientemente grande, de ser o bien disjunta con $P_n$ o bien contener completamente a $P_n$ como subconjunto (so pena de inducir una funci\'on de elecci\'on en un subconjunto infinito de $\{P_n\big|n\in\mathbb N\}$). En otras palabras, debe de existir un $N\subseteq\mathbb N$ y un conjunto finito $F\subseteq X$, que satisface $F\cap P_n=\varnothing$ siempre que $n\in N$, y tal que $A$ es equipotente a $F\cup\bigcup_{n\in N}P_n$.

Esto nos permitir\'a mostrar que el conjunto de cardinales $|A|$ tales que $G$ es $A$-coloreable carece de elemento m\'{\i}nimo. Pues si $c:X\longrightarrow A$ fuera una coloraci\'on adecuada, entonces, como se ha dicho m\'as arriba, podemos suponer sin p\'erdida de generalidad que $A=F\cup\bigcup_{n\in N}P_n$ para alg\'un $F\subseteq X$ finito, y para alg\'un $N\subseteq\mathbb N$ tal que $F\cap P_n=\varnothing$ siempre que $n\in N$. Sea $M=N\setminus\{\min(N)\}$ y note que $G$ tambi\'en es $B$-coloreable, en donde $B=F\cup\bigcup_{n\in M}P_m$ (simplemente redefina $c$ en los $x\in X$ tales que $c(x)\in P_k$, con $k=\min(N)$; esto es posible ya que \'unicamente existen una cantidad finita de tales elementos $x\in X$). Sin embargo, es posible demostrar\footnote{Esto se debe a que todo subconjunto de un conjunto de Russell es {\it Dedekind-finito}, lo que significa que no tiene subconjuntos equipotentes a $\mathbb N$ y, equivalentemente, no admite funciones inyectivas que no sean suprayectivas. El lector interesado puede atacar, como ejercicio, que si $A$ admitiera un subconjunto numerable, entonces ser\'{\i}a posible definir una funci\'on de elecci\'on en alg\'un subconjunto infinito de $\{P_n\big|n\in\mathbb N\}$.} que no puede haber una funci\'on biyectiva $f:A\longrightarrow B$ (pues esta constituir\'{\i}a una funci\'on de $A$ en $A$ inyectiva pero no suprayectiva, ver la nota al pie de p\'agina), por lo cual $|B|\neq |A|$. Como adem\'as $B$ admite una inyecci\'on en $A$ (simplemente t\'omese la funci\'on inclusi\'on), la conclusi\'on es que $|B|<|A|$. Por lo tanto, para cualquier coloraci\'on de $G$, es posible encontrar otra coloraci\'on con una cantidad estrictamente menor de colores, y as\'{\i} el n\'umero $\chi(G)$ no est\'a definido.

\hfill$\Box$

\bigskip
{\bf \centerline{Agradecimientos}}
El primer autor recibi\'o apoyo parcial por parte del proyecto SIP-20221862 del IPN; el segundo autor realiz\'o parte del presente trabajo como becario BEIFI bajo el mismo proyecto. Finalmente, el tercer autor colabor\'o con el presente art\'{\i}culo como parte del programa Delf\'{\i}n de investigaci\'on de verano de 2022.


\bigskip
\hfill\
{\footnotesize
\parbox{5cm}{David J. Fern\'andez-Bret\'on\\
{\it Escuela Superior de F\'{\i}sica y Matem\'aticas},\\
Instituto Polit\'ecnico Nacional,\\
Av. Instituto Polit\'ecnico Nacional s/n Edificio 9, 
Col. San Pedro Zacatenco, Alcald\'{\i}a Gustavo A. Madero, 07738, CDMX, Mexico,\\
{\sf dfernandezb@ipn.mx}}\
{\hfill}\
\parbox{5cm}{Jes\'us A. Flores Hinostrosa\\
{\it Escuela Superior de F\'{\i}sica y Matem\'aticas},\\
Instituto Polit\'ecnico Nacional,\\
Av. Instituto Polit\'ecnico Nacional s/n Edificio 9, 
Col. San Pedro Zacatenco, Alcald\'{\i}a Gustavo A. Madero, 07738, CDMX, Mexico,\\
{\sf jfloresh1501@alumno.ipn.mx}}}\
\hfill

\bigskip
\hfill\
{\footnotesize
\parbox{5cm}{V. Adrián Meza-Campa\\
{\it Facultad de Ciencias Físico-Matemáticas},\\
Universidad Autónoma de Sinaloa.\\
{\it Adscripción Actual: }
Centro de Investigación en Matemáticas, A.C. \\
Jalisco S/N, Col. Valenciana, CP: 36023 Guanajuato, Gto, México.\\
{\sf victor.meza@cimat.mx}}\
{\hfill}\
\parbox{5cm}{L. Gerardo N\'u\~nez Olmedo\\
{\it Escuela Superior de F\'{\i}sica y Matem\'aticas},\\
Instituto Polit\'ecnico Nacional,\\
Av. Instituto Polit\'ecnico Nacional s/n Edificio 9, 
Col. San Pedro Zacatenco, Alcald\'{\i}a Gustavo A. Madero, 07738, CDMX, Mexico,\\
{\sf lnunezo1800@alumno.ipn.mx}}}\
\hfill



\begin{thebibliography}{999}

\bibitem{aharoni:konig}
R. Aharoni,
{\em K\"onig's duality theorem for infinite bipartite graphs.}
J. Lond. Math. Soc. (2)
{\bf 29}~no. 1 (1984), 1--12.

\bibitem{aharoni-matching}
R. Aharoni,
{\em Infinite matching theory.}
Discrete Math.
{\bf 95} (1991), 5--22.

\bibitem{aharoni-berger}
R. Aharoni y E. Berger,
{\em Menger's theorem for infinite graphs.}
Invent. Math.
{\bf 176} (2009), 1--62.

\bibitem{bondy}
Bondy, J.A. y Murty, U.S.R.
{\em Graph Theory with Applications.}
American Elsevier, New York, 1976.

\bibitem{cohen-stcontinuum}
P.J. Cohen,
{\em Set theory and the Continuum Hypothesis.}
W.A. Benjamin, New York/Amsterdam, 1966.

\bibitem{hall-marriage}
P. Hall, 
{\em On representatives of subsets.}
J. London Math. Soc.
{\bf 10} (1935), 26--30.

\bibitem{howard-rubin}
P. Howard y J. Rubin,
{\em Consequences of the Axiom of Choice.}
American Mathematical Society, 1991.

\bibitem{jech-choice}
T. Jech,
{\em The Axiom of Choice.}
North Holland Publishing Company, Amsterdam-London, 1973.

\bibitem{munkres}
J. R. Munkres,
{\em Topology. Second Edition.}
New Jersey, Prentice Hall, 2000.

\end{thebibliography}
\end{document}